\documentclass[mnsc,nonblindrev,copyedit]{informs3}

\usepackage{graphicx}
\usepackage{multirow}
\usepackage{float}
\usepackage{setspace,graphics,booktabs,color,enumerate,subfigure}
\usepackage{amssymb,amsmath}
\newtheorem{proposition}{Proposition}[section]
\usepackage{mathrsfs}
\usepackage{url}
\definecolor{DarkBlue}{rgb}{0,0.08,0.45}
\usepackage[backref = false, bookmarks, colorlinks = true, plainpages = false, citecolor = DarkBlue , urlcolor = DarkBlue, filecolor = DarkBlue, linkcolor = DarkBlue]{hyperref}

\usepackage{enumerate}
\usepackage{cases}
\usepackage{epstopdf}
\usepackage[linesnumbered, ruled]{algorithm2e}
\usepackage{threeparttable}
\usepackage{IEEEtrantools}

\SetKwRepeat{Do}{do}{while}%

\usepackage{natbib}
 \def\bibfont{\small}%
\theoremstyle{TH}%

\theoremstyle{EX}
\newtheorem{rem}{Remark}
\newtheorem{exmp}{Example}

\newtheorem{defn}{Definition}

\ECRepeatTheorems

\EquationsNumberedThrough    % Default: (1), (2), ...
\MANUSCRIPTNO{}

\begin{document}

\begin{center}
{\LARGE \bf{ Target-based Distributionally Robust Minimum Spanning Tree Problem}}
\end{center}

% \footnotetext{University of Chicago}
% \footnotetext{Nanjing University}
\AUTHOR{%
  Yang Xu\thanks{Univeristy of Chicago, yxu6@uchicago.edu} 
  Lianmin Zhang\thanks{Nanjing University}}
  
{\par\noindent\bf Abstract:}
Due to its broad applications in practice, the minimum spanning tree problem and its all kinds of 
variations have been studied extensively during the last decades, for which a host of efficient exact 
and heuristic algorithms have been proposed. Meanwhile, motivated by realistic applications, the minimum spanning tree problem in stochastic network has attracted considerable attention of researchers, with respect to which stochastic and robust spanning tree models and related algorithms have been continuingly developed. However, all of them would be either too restricted by the types of the edge weight random variables or computationally intractable, especially in large-scale networks. In this paper, we introduce a target-based distributionally robust optimization framework to solve the minimum spanning tree problem in stochastic graphs where the probability distribution function of the edge weight is unknown but some statistical information could be utilized to prevent the optimal solution from being too conservative. We propose two exact algorithms to solve it, based on Benders decomposition framework and a modified classical greedy algorithm of MST problem (Prim algorithm),respectively. Compared with the NP-hard stochastic and robust spanning tree problems,The proposed target-based distributionally robust minimum spanning tree problem enjoys more satisfactory algorithmic aspect and robustness, when faced with uncertainty in input data.
%  \end{abstract}
    
\begin{flushleft}{\bf Key words:}
Minimum Spanning tree ; Distributionally robust optimization ; stochastic graph.
\end{flushleft}

\section{Introduction}

The minimum spanning tree (MST) problem is one of the most well-studied network flow problems. Given an undirected network $G$ with real-valued 
edge weight, a minimum spanning tree (MST) is a spanning tree with sum of edge weights minimum among all spanning trees of $G$.
Extensive exact and heuristic algorithms, such as branch-and-cut, greedy algorithm and tabu search algorithm, are commonly used for optimizing
MST solutions and some of its variants, when edge weights are deterministic. These problems and their solutions have proven to be 
essential ingredients for addressing many real-world problems in applications as diverse as transportation, computer networking and so on.

In many of these applications, the presence of uncertainty regarding edge weights in the networks is a critical issue to consider explicitly
if one hopes to provide solutions of practical values to the end users. In the
past, the main work of research on MST problems on stochastic graph, despite excellent, would be unsuitable for the practical use
because of the intractable computational cost or failure to deal with stochastic graphs with edge
weight random variables subject to continuous or even unknown probability distribution.
Recently, regarding to the presence of different kinds of uncertainty in networks,
some studies proposed the target driven theory and the distributionally robust optimization to deal with them.
The target driven theory would bypass the use of utilities in the classical stochastic optimization to reduce subjectivity
and also fit more for reality. The distributionally robust optimization, then,   
would avoid the solution being too conservative through involving some statistical information when compared with the original 
robust optimization framework and enjoy better computational aspect and 
robustness than the stochastic optimization when faced with uncertainties, especially with the use of certain objective functions 
having satisfying properties such as monotonicity and convexity as index.
Thus, in this paper, we introduce the target-based distributionally robust optimization framework
to deal with the MST problems on stochastic graphs, leading to our target-based distributionally robust minimum spanning tree (TDRMST)
model.

We consider the minimum spanning tree problem with edge weight under uncertainties. In this case, the exact 
distribution of edge weight is unknown and some statistical information, such as bounded support and moments, is available. The objective
is to minimize the \textit{Requirement Violation(RV) index} which in this case represents the risk of outstriping the given target set 
for the whole weight under certain ambiguity set. The main contributions of the paper are as below.

\begin{itemize}
	\item We discuss the extensions of the distributional robust optimization framework to the case of MST problem, giving out a convex mixed integer programming formulation of the model and some analysis of the properties of the \textit{RV index}.
	\item Based on the properties, we propose efficient algorithm based on standard MST algorithm called
	Prim algorithm, which could outperform both the original Benders Decomposition framework and the bisection method to solve the 
	target-based distributionally robust minimum spanning tree (TDRMST) problem.
\end{itemize}

The rest of our paper is organized as follows. In Section 2, we summarize the related literature. In Section 3, we introduce the target-based decision theory and the \textit{RV Index} to quantify
the risk of violating the target. In Section 4, we present an optimization framework for the problem and propose two efficient solution procedures to solve the problem exactly. In Section 5, we perform 
various numerical experiments to illustrate the benefits of our model. Finally, we draw conclusions in Section 6.

\section{Literature Review}

In this section, we summarize the related literature in regard to both the deterministic MST problem with its variants and its stochastic versions
which could be divided into stochastic MST problems and robust MST problems. Also, we introduce the target driven optimization and distributionally
robust optimization for the later proposition of our TDRMST model.

\subsection{Minimum spanning tree problem}

A Minimum spanning tree of a weighted, undirected graph $G$ is a spanning tree
of $G$ whose edges sum to minimum weight. With so many applications in practice such as the design of communication network, electric power system and so on, the minimum spanning tree problems and its variants
have been studied extensively over many decades (\cite{akbari2012learning}).  
	
In general, MST problems could be subdivided into deterministic MST problems and stochastic MST problems, based on whether the edge weight is 
assumed to be fixed or random variable. For a standard deterministic MST problem, greedy algorithms like \cite{10.2307/2033241} and 
\cite{prim1957shortest} algorithm would be used to build an optimal spanning tree. Moreover, its all kinds of variants like MST problems with additional constraints or other combinatorial optimization problems based on the context of spanning tree, are continuingly developed and corresponding algorithms are proposed.

\subsubsection{Variants of MST problem}

In several practical contexts like communication systems, a problem of interest
would be to find a minimal spanning tree subject to additional constraints, like
$\mathop{\sum}\limits_{j:e_j \in E_T} c_j \le L,$

As is shown in \cite{aggarwal1982minimal}, this kind of problems belongs to the class of weakly NP-hard problems. 
Representative constrained minimum spanning tree problems studied in the literature are summarized in Table 1.1:

\begin{table}[ht!]
\begin{center}
	\caption{Summary of the Related Literature}
	\resizebox{\textwidth}{!}{
	\begin{tabular}{c|c|c}
		\hline
		\textbf{Author(s)} & \textbf{Problem} & \textbf{Approach}\\
		\hline
		\cite{gruber2006neighbourhood} & Bounded Diameter MST(BDMST) & Heuristic(VNS/EA/ACO)\\
		\cite{bui2006ant} & Degree constrained MST(DCMST) & Heuristic(ACO)\\
		\cite{oncan2007design} & Capacitated MST(CMST) & Heuristic\\
		\cite{oncan2008tabu} & Generalized MST(GMST) & Heuristic(TS)\\
		\cite{parsa1998iterative} & Delay-Constrained MST & Heuristic\\
		\cite{gouveia2011modeling} & Hop-Constrained MST(HMST) & Branch-and-cut\\
		\hline
	\end{tabular}}
\end{center}
\end{table}

To be specific, the eccentricity of a node $v$ is the maximum number of edges on a path
from $v$ to any other node in the tree $T$. The diameter of $T$ is the maximum eccentricity of 
all nodes. Thus BDMST problem is a combinatorial optimisation problem searching for a spanning tree rooted 
at an unknown center having its height restricted to half of the diameter and \cite{gruber2006neighbourhood}
explores various heuristic algorithms such as variable neighbourhood search (VNS), an evolutionary
algorithm (EA) and an ant colony optimisation (ACO) to search for solutions. DCMST problem tends to search for
a minimum-cost spanning tree such that no vertex in the tree exceeds a given degree constraint as its 
name implies. CMST problem endeavors to find a spanning tree of minimum cost so that the total demand of the vertices
in each subtree rooted at the central vertex would not exceed the capacity limitation. \cite{oncan2007design} 
introduces fuzzy input data into the classical CMST problem and designs a fast approximate reasoning algorithm based on
heuristic algorithms and fuzzy logic rules to solve it. GMST problem consists of designing a minimum-cost tree spanning all clusters, 
into which the vertex set of a graph is partitioned. Aimed at the GMST problem, \cite{oncan2008tabu} develop a tabu search (TS) algorithm to
solve it. Moreover, Delay-Constrained MST problem is to
construct minimum-cost trees with delay constraints and HMST problem is to find a minimum-cost tree such that the unique path
from a specified root note to any other node has restricted number of edges(hops). Given that the constrained MST problems 
are NP-hard problem, heuristic methods have been extensively used by researchers for solving those problems as is shown in the
Table 1.1.

Besides constrained MST problems, other variants of MST problems have caught more and more attentions. 
\cite{sokkalingam1999solving} first studies the inverse spanning tree problem and formulates it 
as the dual of an assignment problem on a bipartite network. Taking advantage of the formulation, they
implement the successive shortest path algorithm which could run in $O(n^3)$ time. 
\cite{hochbaum2003efficient} proposes more efficient algorithms for the inverse spanning tree problem, whose run time
could be $O(nm\log^2n)$ for any convex inverse spanning-tree problem, where $n$ denotes the number of nodes and $m$ represents the number of edges.
\cite{doi:10.1287/ijoc.1040.0123} explores the Multilevel Capacitated MST (MLCMST) problem, a generation of the well-known CMST problem, which allows 
for multiple facility types in the design of the network. They develop flow-based mixed integer programming formulations 
to find tight lower bounds and develop heuristic procedures for the MLCMST problem. \cite{sourd2008multiobjective}
summarizes the related work of the Bi-objective Minimum Spanning Tree (BOST) problem, which is to find one spanning 
tree for each Pareto point in the objective space. This paper implement a multi-objective Branch-and-Bound framework
to the BOST problem and shows the efficiency of the approach. More recently, \cite{wei2021integer} considers a two-player 
interdiction problem staged over a graph where the attacker's objective is to minimize the cost of removing edges
from the graph so that the weight of a minimum spanning tree in the residual graph would be increased up to a predefined
level $r$. The authors provide a detailed study of the problem's solution space, present multiple integer programming
formulations and a polyhedral analysis of the convex hull of feasible solutions. \cite{paul2020budgeted} considers versions of
the prize-collecting MST problems and presents 2-approximation algorithm based on a parameterized primal-dual approach.

\subsubsection{Stochastic spanning tree problem}
Despite the extensive researches on the deterministic MST problems and related algorithms, in actual situations like the construction of a communication network which connects some cities, the edge weights or costs could vary with time, leading to the significance of considering the stochastic version of MST problem where edge weights are not constant but random variables. On the whole, the stochastic spanning tree problem has two main stochastic versions.

On one hand, it's about the chance constraint models. \cite{ishii1981stochastic} firstly generalizes the minimal spanning tree problem toward a 
stochastic version, considering a stochastic spanning tree problem in which edge costs are not
constant but random variables, and its objective is to find an optimal spanning tree satisfying a 
certain chance constraint, i.e.
\begin{equation}
\begin{aligned}
	\min &\; f, \\
	\text{s.t.} &\; Prob\{\sum_{j=1}^{m} c_jx_j \le f\} \ge \alpha,\\
	&\; x_j \in \{0,1\}, \quad X: spanning \ tree, \\
\end{aligned}	
\end{equation}
\cite{ishii1981stochastic} assumes $c_j$ to be distributed according to the normal 
distribution $N(\mu_j, \sigma_j^2)$ and be mutually independent. They also transform 
the problem into its deterministic equivalent problem and propose an auxiliary problem, based on 
which they propose a polynomial order algorithm to attain exact solution.\cite{ishii1983stochastic}
extend the similar definition to the bottleneck spanning tree problem. Based on these work,
\cite{mohd1994interval} and \cite{ishii1995confidence} further make analysis and improvements.
\cite{mohd1994interval} introduces the method called an interval elimination to solve the 
stochastic spanning tree problem while \cite{ishii1995confidence} introduces the confidence regional
method to deal with the uncertainty in the unknown parameters of underlying probability distribution.
More recently, \cite{shen2015chance} considers a balance-constrained stochastic bottleneck spanning tree (BCSBST) problem.
The paper formulates the problem as a mixed-integer nonlinear program, develops two 
mixed-integer linear programming approximation and proposes a bisection algorithm to approximate optimal solutions in polynomial time.

On the other hand, the stochastic minimum spanning tree (SMST) is defined as a stochastic spanning tree with the minimum
expected weight, i.e. stochastic spanning tree $\tau^* \in T$ is the stochastic MST if and only if 
$\overline{w}_{\tau^*} = \min_{\forall \tau_{i}\in T}\{\overline{w}_{\tau_i}\}$ where $\overline{w}_{\tau_i}$
denotes the expected weight of spanning tree $\tau_i$.

With regard to this kind of stochastic version, several scholars concentrated their attention on
establishing bounds on the distribution of $W$ and making asymptotic analysis.
\cite{jain1988approximations} firstly proposed a method to obtain bounds, which is shown to be much tighter than the naive bound 
obtained by computing the MST length of the deterministic graph with the respective means as arc lengths, and approximations
for the MST length on an undirected graph whose arc lengths are independently distributed random variables. They also
analyze the asymptotic properties of their approximations. \cite{hutson2005bounding} comprehensively summarized the related work
in the analysis of obtaining tighter bounds and distribution of the MST lengths on graphs under different conditions and make some improvement.
It also deserves to be mentioned that \cite{alexopoulos2000state} investigated state space partitioning technique, which is considered in 
\cite{doulliez1972transportation}, to compute and bound specific values of the minimum spanning tree distribution in graphs
with independent, but not necessarily identically distributed, discrete edge weight random variables.
Based on this partitioning technique, \cite{hutson2005bounding} employed a specific heuristic approach
to obtain nonintersecting sets from fundamental cutsets and cycles to derive better bounds on $E[W]$.
More recently, \cite{frieze2021randomly} studied the minimum spanning tree problem with an additional constraint
that the optimal spanning tree $T$ must satisfy $C(T)\le c_0$ on the complete graph with independent uniform edge weight random variables $U^s$.
They mainly establish the asymptotic value of the optimum weight for a range of $c_0$ and $s$ through the consideration of a dual problem.

Several authors, otherwise, have also devoted themselves to obtaining better solutions where the edge weights are determined by different 
kinds of random variables.
To solve the MST problem in networks where the edge weight can assume a finite number of distinct values,
\cite{hutson2006minimum} considered several approaches such as repeated prim (RP) method, cut-set(CT) method, cycle tracing method (CTM) and 
multiple edge(ME) sensitivity method, among which the best method has a worst case running time of $O(mN)$, 
where $m$ is the number of edges and $N$ is the number of states.
For the situation when the edges are continuous random variables with unknown probability distribution function,
\cite{he2008model} proposed a hybrid intelligent algorithm as a combination of the genetic algorithm and stochastic simulation,
taking advantage of the Prufer encoding scheme to represent all trees to code the corresponding spanning tree
for the genetic representation. \cite{akbari2012learning} and \cite{torkestani2011learning} employed a learning 
automata-based heuristic algorithm which could decrease the number of samples
that need to be taken from the graph by reducing the rate of unnecessary samples.
Furthermore, \cite{de2005evolutionary} introduced the Fuzzy Set Theory into the MST problem to deal with the uncertainty and examines
the MST problem with fuzzy parameters, in case of which they propose both exact and genetic algorithms to solve it.

\subsubsection{Robust spanning tree problem}

The same as stochastic optimization, robust optimization also aims at dealing with uncertainty in input data. 
Nevertheless, stochastic optimization, to deal with the unknown parameter, is in need of the underlying probability distribution model
or its assumption, which could be hardly possible in reality. 
Thus, robust optimization, which could bypass the precondition and only need some statistical information, gradually draws scholars' attention.
Concretely speaking, robust optimization is developed to hedge
against the worst possible scenario according to the given degree of conservative and some information
of the data with unknown distribution, such as all the cases of the discrete parameters or 
the interval data.

Motivated by telecommunications applications, authors investigate the minimum spanning tree
problem on networks whose edge weights or costs are interval numbers, and define robust spanning tree problem
to hedge against the worst case contingency, under the robust deviation framework. 

Generally, the goal of the robust spanning tree problem with interval data (RSTID) is to find a spanning tree which
minimizes the maximum deviation of its cost from the costs of the optimal spanning trees obtained for all possible
realizations of the edge costs/weights within the given intervals, called a relative robust spanning tree as denoted 
in \cite{yaman2001robust}. That is $T^r \in \arg\min_{T:spanning \  tree}d_T$, where 
$d_T = \max\{c_T^s-c_{T*(s)}^s\}$ is the robust deviation for spanning tree $T$ for scenarios $s$. 

\cite{yaman2001robust} presented a mixed integer programming formulation for the problem and define some useful optimality
concepts like weak edges and strong edges. Then they present characterizations for these entities, which are used to preprocess
a given network with interval data prior to the solution of the mixed integer programming problem. They present methods to 
identify all weak and strong edges in the graph in polynomial time, leading to the exclusion of certain edges to be in the optimal spanning tree. 
\cite{aron2004complexity} later studied the complexity of the RSTID problem and proves it to be NP-complete. \cite{montemanni2005branch} 
and \cite{montemanni2006benders} proposed a Branch-and-bound and Benders Decomposition approach respectively to solve the problem more efficiently.
\cite{salazar2007robust} gave some characterizations of strongly strong edges and non-weak edges leading to recognition algorithms.

Moreover, researchers have developed other models, which take advantage of other forms of the objective function to deal with the uncertainty in MST problem on stochastic networks.
Related to the RSTID problem, \cite{chen2009polynomial} proposed and studied a new model for the spanning tree problem with interval data, called
the Minimum Risk Spanning Tree (MRST) problem which is to establish a spanning tree $T$ of total cost not more than a given constant so that the risk sum
over the links in $T$ is minimized, where the risk is $\frac{\overline{c_e}-c_e}{\overline{c_e}-\underline{c_e}}$ for $c_e \in [\underline{c_e},\overline{c_e}]$.
As is shown in \cite{chen2009polynomial}, the MRST model could be solved in polynomial time and more satisfactory in algorithmic aspect than RSTID model.Plus, \cite{li2011maximizing} study the stochastic versions of a broad class of combinatorial problems, including the minimum spanning tree problem,
where the weights of the elements in the input dataset are uncertain. Their general objective is to maximize the expected utility of the solution
for some given utility function, rather than the expected weight and obtain a polynomial time approximation algorithm with additive error $\epsilon$
for any $\epsilon > 0$ when the conditions are not bad.

To put it in a nutshell, the past work, despite excellent, would be unsuitable for realistic applications when faced with stochastic networks. 
Stochastic spanning tree problems need to assume the probability distribution function which could be inaccessible in reality. 
Robust spanning tree problems are too conservative because of failure to involve more statistical information. 
Both of them could also be confronted with intractable computational cost, especially in large-scale networks.

\subsection{Target driven optimization}
The classical approaches to deal with uncertainty would either choose to minimize the expected profits over a multistage planning horizon 
(\cite{hutson2006minimum},\cite{he2008model},\cite{akbari2012learning},\cite{torkestani2011learning},\cite{de2005evolutionary}), 
which is suitable for risk neutral decision makers, or address risk through optimizing over a mean risk objective 
(\cite{chen2009polynomial}), or an expected utility (\cite{li2011maximizing}).
The former requires numerous number of repetition of the decision under identical conditions, which could be impossible at the most time. 
The latter would be too subjective with respect to the articulation of the decision makers' utility function, 
making it hard to ascertain in practice. 
Recently, researchers in decision theory introduces the definition of an aspiration level, 
target or goal to bypass the use of utility functions (\cite{chen2009goal}). 
The objective of target driven optimization is generally to maximize the probability of reaching the target, 
which has been playing a significant role in daily decision making since actually, decision makers are more 
concerned about the prospect of failure to reach some target rate of return rather than the risk itself (\cite{mao1970survey}).

\subsection{Distributionally robust optimization}

In spite of the ability to deal with data with unknown distribution, robust optimization usually suffer from 
too conservative solutions and computational cost. In this case,the framework of distributionally robust optimization 
has drawn more and more attention, especially with the proposition of the solution 
frameworks which could incorporate some statistical information, extending the limited input data like 
all the cases of the discrete parameters or the interval data, into the models to find a solution with a high
level of robustness while not overly conservative (\cite{bertsimas2011theory}).
\cite{rahimian2019distributionally} gave a more general review.

Nowadays, the framework of the distributionally robust optimization has been extensively applied and improved
in the vehicle routing problems and other optimization problems. With respect to a routing problem with soft time windows 
where exact probability distributions of travel times are not known for a single 
uncapacitated vehicle, \cite{jaillet2016routing} proposed a performance index, termed as the \textit{Requirements Violation (RV) Index} 
that represents the risk of violating the time window restrictions and a general solution framework. They also discussed several 
special cases and provided computational results showing that the proposed performance index could produce solutions generally 
superior to other approaches, including stochastic programming solved by sampling techniques. Moreover, 
under such distributionally robust framework,
\cite{long2014distributionally} studied the discrete optimization problem by optimizing the Entropic Value-at-Risk and 
proposed an efficient approximation algorithm.
\cite{cui2021inventory} researched on an uncertain inventory routing problem 
and provided exact algorithms to solve the problem under the Service Violation Index.
However, the distributionally robust optimization, as far as we know, has not been introduced into the MST problems.

\section{Model Formulation}

In this section, we build our model to deal with the minimum spanning tree problems in stochastic networks with target driven policy. 
Firstly, we give out the following notations.

%\begin{note}
Let $G=(V,E)$ be an undirected connected graph with vertexes $V=\{v_1,v_2,...,v_n\}$ and edges $E=\{e_1,e_2,...,e_m\}$. 
Each edge $e_j$ is endowed with a weight $w_j$, $j=1,\cdots,,m$.
We take boldface characters to represent vectors, e.g., $\boldsymbol{w}=(w_j)_{j\in[m]}$, where $[m]=\{1,\dots,m\}$. 
Plus, we denote random variables by characters with tilde sign $\tilde{\cdot}$, e.g., 
$\tilde{w_j}$ denotes the edge weight random variable. With $\mathscr{V}$ as the space of real-valued random variables, we define 
the sample space, the associated sigma algebra and the true distribution as $(\Omega, \mathbb{F}, \mathbb{P})$, respectively. 
In our case where we have no access to the full information for $\mathbb{P}$, 
we assume that $\mathbb{P} \in \mathbb{F}$ and call $\mathbb{F}$ the uncertainty set. We denote the 
expectation of $\tilde{\boldsymbol{w}}$ subject to the probability distribution $\mathbb{P}$ as $E_{\mathbb{P}}[\tilde{\boldsymbol{w}}]$.
%\end{note}

A connected acyclic spanning subgraph $T=(V,E_T)$ of $G$ is called a \textit{spanning tree} of $G$ 
if it satisfies $E_T \subseteq E$ and $|E_T|=n-1$.  More specifically,

\begin{defn}
	A subgraph $T=(V,E_T)$ of $G$ is called a \textit{spanning tree} of $G$ if it satisfies 
	\begin{itemize}
		\item[1)] Subgraph $T=(V,E_T)$ is connected.
		\item[2)] $T=(V,E_T)$ has the same vertex-set as $G$.
		\item[3)] $|E_T|=n-1$.
	\end{itemize}
\end{defn}

The weight of a spanning tree is given by $w(T)=\mathop{\sum}\limits_{j:e_j \in E_T} w_j$.
Denote the set of all spanning trees of $G$ by $\mathscr{T}(G)$. Then the classical Minimum Spanning Tree (MST) Problem 
can be stated as $\mathop{\min}\limits_{T\in\mathscr{T}(G)}w(T)$.
We choose the following network design type of formulation (\cite{magnanti1995optimal},\cite{wei2021integer}) that is commonly used to characterize the spanning trees in a graph $G$:
\begin{equation}
	\begin{aligned}
		\min &\; \sum_{e \in E}  w_e y_e, \\
		\text{s.t.} &\; \sum_{j:{i,j} \in E}f_{ij}^l - \sum_{j:{i,j} \in E}f_{ji}^l  =b_i^l, &\, \forall i \in V,l \in V\setminus{k}, \\
		&\; y_e \ge f_{ij}^l+f_{ji}^l, &\, \forall e =\{i,j\}\in E,l \in V\setminus{k},\\
		&\; \sum_{e \in E}y_e  =n-1,\\
		&\; f_{ij}^l,f_{ji}^l  \ge 0, &\, \forall e=\{i,j\}]\in E,l \in V\setminus{k}, \\
		&\; y_e  \ge 0,  &\, \forall e\in E
	\end{aligned}	
\end{equation}

where $\boldsymbol{y} = {[y_e]}_{e \in E}$ represents the incidence vector of the spanning tree to be selected. 
In this setting, a vertex $k \in V$ is arbitrarily chosen to be the root of the spanning tree 
and a commodity $l$ is defined for every vertex $l \in V\setminus \{k\}$. 
A parameter $b_i^{l}$, which is used to model the requirement that each unit of commodity $l$ generated from the center $k$ 
must be sent to vertex $l$, selects the values of 1, -1, and 0 when $i = k, i = l, $ and $i \in V\setminus \{k, l\},$ respectively.
$f_{ij}^l$, called the flow decision variable, denotes the flow of the commodity $l$ from $i$ to $j$, 
which is commonly used in dealing with such problems. The first few constraints describe the balance flow constraints 
so that it would be connected and the restriction of the number of edges to be $n-1$ is added to make it a spanning tree. 

\begin{rem}
Note that although it would be more complicated when we introduce the flow decision variable, the formulation above indicates the fact that $\boldsymbol{y}$ would always be in $\{0,1\}$ without needing to enforce integrality constraints \cite{magnanti1995optimal}.Plus, as is said in \cite{wei2021integer}, the set of extreme points of this formulation corresponds to solutions $(\boldsymbol{y},\boldsymbol{f})$. This formulation would also be beneficial when we extend the work from the standard MST problem to its other variants, which would generally use this formulation to represent the spanning tree configuration. 
\end{rem}

%origin

However, in many practical settings the weights $w_j(j=1,\cdots,,m)$ are random rather than deterministic, following some unknown probability distribution. The traditional models of minimum spanning tree are invalid in this scenario. To fulfill this gap, we formulate and solve a target-based distributionally robust minimum spanning tree model in this study. Let the random weights $\tilde{\boldsymbol{w}}$ follow a distribution $\mathbb{P}$, which belongs to the uncertainty set
$\mathbb{F} =\left\{
\mathbb{P} \;\left|\; E_{\mathbb{P}}[\tilde{\boldsymbol{w}}] = \boldsymbol{\mu},
\underline{\boldsymbol{w}} \leq \tilde{\boldsymbol{w}} \leq \overline{\boldsymbol{w}}
\right.
\right\}$,
where $\boldsymbol{\mu},\underline{\boldsymbol{w}},\overline{\boldsymbol{w}} \in\mathbb{R}_+^m$. For any spanning tree $T\in\mathscr{T}(G)$, 
its weight is also a random variable $\tilde{w}(T)= \mathop{\sum}\limits_{j}\tilde{w}_j y_{e_j\in E_T}$. 
Denote $\boldsymbol{y} \in \mathscr{Y} = \{\boldsymbol{y}_T = (y_{e_j\in E_T})_{j=1}^m | T\in\mathscr{T}(G)\}$. 
Then $\tilde{w}(\boldsymbol{T})=\tilde{w}(\boldsymbol{y}_T)= \tilde{\boldsymbol{w}}'\boldsymbol{y}_T$

i.e.

$$	
	\mathscr{Y}=\Biggl\{\boldsymbol{y} \in {\{0,1\}^{m}}\left|
	\begin{array}{l}
	\sum_{j:{i,j} \in E}f_{ij}^l - \sum_{j:{i,j} \in E}f_{ji}^l  =b_i^l, \forall i \in V,l \in V\setminus{k}, \\
    y_e \ge f_{ij}^l+f_{ji}^l,  \forall e =\{i,j\}\in E,l \in V\setminus{k}, \\
    \sum_{e \in E}y_e  =n-1,\\
    f_{ij}^l,f_{ji}^l  \ge 0, \forall e=\{i,j\}\in E,l \in V\setminus{k}.
	\end{array} \right.\Biggr\}
$$

\subsection{Requirements Violation Index}

Target-based decision argues that the main goal of most firms is to attain the target(s) rather than optimizing objectives (\cite{Hall2015}). 
Given a target $\tau$, the difference $\tau -\tilde{\boldsymbol{w}}'\boldsymbol{y}$ means the gap between the exact cost and the target, 
which can be explained as the target premium (\citet{Chen2009}). Surely, we expect that 
$\mathbb{P}(\tau \ge\tilde{\boldsymbol{w}}'\boldsymbol{y}) = 1$, which implies that the target can be achieved almost surely. 
However, sometimes the situation could be impossible. Thus we can intoduce some performance index, called the \textit{RV Index}, 
whose definition is based on the idea of certainty equivalent. 

\begin{defn}
	\textbf{\emph{Requirements Violation (RV) Index:}} Given an uncertain attribute $\tilde{t}$ and its lower and upper limits, $\underline{\tau},\overline{\tau}$, the RV Index $\rho_{\underline{\tau},\overline{\tau}}(\tilde{t}):\mathscr{V}\to [0, +\infty]$ is defined as follows {\cite{jaillet2016routing}}:
	$$
	\rho_{\underline{\tau},\overline{\tau}}(\tilde{t})=\inf \left\{ \alpha : C_{\alpha}(\tilde{t}) \leq \overline{\tau},C_{\alpha}(-\tilde{t}) \leq -\underline{\tau}, \alpha \ge 0 \right\},
	$$
	or $+\infty$ if no such $\alpha$ exists, where $C_{\alpha}(\tilde{t})$ is the worst-case certainty equivalent under exponential disutility defined as 
	$$C_{\alpha} (\tilde{t}) = \alpha ln \mathbb{E}_{\mathbb{P}}\left[exp\left(\frac{\tilde{t}}{\alpha}\right)\right], \textit{if} \  \alpha >0,$$
	or 
	$$C_{\alpha}(\tilde{t})=\lim_{\gamma \downarrow 0}C_{\gamma}(\tilde{t}), \textit{if} \  \alpha =0.
	\footnote{Here, it corresponds to the worst-case deterministic value which an individual under Constant Absolute Risk Aversion(CARA) of risk tolerance parameter $\alpha\ge0$ would equally prefer over an uncertain attribute $\tilde{t}$. We strongly recommend interested reader to reach \cite{chen2021robust},\cite{jaillet2016routing} for more detailed analysis with respect to the index.}$$
\end{defn}

\begin{rem}
	When $\tilde{t}$ is constant, we have $C_{\alpha}(\tilde{t})=\textit{constant}$; And if the probability distribution of random variable $\tilde{t}$ is accessible, we could obtain the value of $C_{\alpha}(\tilde{t})$. For example, if $\tilde{t} \sim N(\mu, \sigma^2)$, we have that 
   $$C_{\alpha}(\tilde{t})=\alpha \ln \mathbb{E}_{\mathbb{P}}(\exp(\frac{\tilde{t}}{\alpha}))=
   \alpha \ln(\exp(\frac{\mu}{\alpha}+\frac{\sigma^2}{2\alpha^2})
   =\mu + \frac{1}{2\alpha}\sigma^2.$$
\end{rem}

The RV Index enjoys the \textbf{Full satisfaction} property, which means that: $\rho_{\underline{\tau},\overline{\tau}}(\tilde{t}) = 0$ if and only if $\mathbb{P}(\tilde{t}\in [\underline{\tau},\overline{\tau}])=1$ for all $\mathbb{P}\in \mathbb{F}$. It implies that it's most preferred when $\rho =0$ and when $\rho =+\infty$, the uncertain attribute is strongly suggested to be move out from later consideration. 

And the RV Index could be easily extended to collective RV Index by several objective functions on it. For example, $\sum_{i\in \mathcal{I}}\rho_{\underline{\tau_i},\overline{\tau_i}}(\tilde{t_i})$ or $\max_{i\in \mathcal{I}}\rho_{\underline{\tau_i},\overline{\tau_i}}(\tilde{t_i})$ could both preserve the satisfactory properties of RV Index while the prior would be apparently easier. The collective RV Index would definitely be profitable for the extension of our work to other MST variant problems with additional constraints, where we could be challenged with the jointly probability of a set of attributes in realizing targets. Also, we could enforce weights to the lower and upper bound requirements in order to differentiate the significance of them, leading to the following formulation: 
$$\rho_{\underline{\tau},\overline{\tau}}(\tilde{t})=\inf \left\{ \alpha : C_{w_1\alpha}(\tilde{t}) \leq \overline{\tau},C_{w_2\alpha}(-\tilde{t}) \leq -\underline{\tau}, \alpha \ge 0 \right\}.$$

In our context, we take its special case where we only enforce one target, i.e. the upper bound rather than both lower and upper bounds and also $\tilde{t} = \tilde{\boldsymbol{w}}'\boldsymbol{y}$. Thus, the deterministic quantity $C_{\alpha} (\tilde{\boldsymbol{w}}'\boldsymbol{y})$ is the certainty equivalent of random travel cost $\tilde{w}(T)$ and $\alpha > 0$ is a risk tolerance parameter associated with failure to realize the target. 

For a given target $\tau$, the quality of the random travel cost $\tilde{w}(T)$ will then be defined as the smallest risk tolerance $\alpha$ allowable 
so that the certainty equivalent of travel cost does not exceed the target $\tau$; i.e.,

$$
\rho_{\tau}(\tilde{\boldsymbol{w}}'\boldsymbol{y})=\inf \left\{ \alpha : \sup_{\mathbb{P} \in \mathbb{F}} \alpha ln\mathbb{E}_{\mathbb{P}}\left[exp\left(\frac{\tilde{\boldsymbol{w}}'\boldsymbol{y}}{\alpha}\right)\right] \leq \tau, \alpha \ge 0 \right\}.
$$

Finding a minimum spanning tree for which the travel cost $\tilde{w}(T)$ gives the smallest policy can thus be obtained by solving the following optimization problem: 
$$\mathop{\min}\limits_{\boldsymbol{y} \in \mathscr{Y}}\rho_{\tau}(\tilde{\boldsymbol{w}}'\boldsymbol{y})$$
or

%origin

\begin{equation}\label{original_problem}
	\begin{aligned}
		\inf &\; \alpha \\
		\text{s.t.} &\; h(\alpha, \boldsymbol{y}) \leq \tau,  \\
		&\; \alpha \ge 0,  \\
		&\; \boldsymbol{y} \in \mathscr{Y}.
	\end{aligned}	
\end{equation}

where
$$h(\alpha, \boldsymbol{y})=C_{\alpha} (\tilde{\boldsymbol{w}}'\boldsymbol{y}_T)=\sup_{\mathbb{P} \in \mathbb{F}} \alpha ln\mathbb{E}_{\mathbb{P}}\left[exp\left(\frac{\tilde{\boldsymbol{w}}'\boldsymbol{y}_T}{\alpha}\right)\right].$$

Different from the chance constrained model and robust formulation stated before, our target-based distributionally robust MST(TDRMST) model would be better in the following ways. Firstly, the former only ensures the probability of violation at certain given level, failing to take into consideration the magnitude, while the RV Index could measure the violation probability at any level. When compared with the latter, the TDRMST formulation could involve several descriptive information so as to prevent the optimal solution from too conservative. And for both of them, the TDRMST formulation enjoys more satisfactory computational aspect because of the following analysis and proposed efficient algorithms.

\section{Solution Procedure}

In this section, we develop the optimization framework for the model and then give out efficient solution procedures.

\subsection{Robust-Optimization(RO) algorithm}

Although solving the problem above would be challenging, if the vector $\overline{\boldsymbol{y}} \in \mathscr{Y}$ is known, the corresponding objective function, denoted by $f^r (\overline{\boldsymbol{y}})$ can be computed by the following 
convex problem, denoted as $SP(\bar{\boldsymbol{y}})$
\begin{equation}
	\begin{aligned}
		f^r (\overline{\boldsymbol{y}}) = \inf &\; \alpha \\
	    \text{s.t.} &\; h(\alpha, \overline{\boldsymbol{y}}) \leq \tau,  \\
	    &\; \alpha \ge 0.	
	\end{aligned}
\end{equation}

Thus, with the following useful conclusions, we give out an algorithm based on the approximation techniques and the Benders decomposition framework to solve the problem.

\begin{proposition}
    For any $\boldsymbol{y} \in \mathscr{Y}$, we have 
    $$f(\boldsymbol{y}) = \sup_{\boldsymbol{s} \in \mathscr{Y}}\{f(\boldsymbol{s})+d^f_{\boldsymbol{s}}(\boldsymbol{s})'(\boldsymbol{y}-\boldsymbol{s})\},$$
    where $d^f_{\boldsymbol{s}}(\boldsymbol{s})$ is the vector of subgradient of $f(\boldsymbol{s})$ with respect to $\boldsymbol{s}$.\footnote{The detailed proof could be found in Proposition 6 \cite{jaillet2016routing}.}
\end{proposition}

This proposition gives a piece-wise linear approximation of convex $f(\boldsymbol{s})$ and makes it possible 
for us to use Benders decomposition method to tackle the original problem.

\begin{proposition}
    The original problem(\ref{original_problem}), target-based distributionally robust MST problem, could be formulated as below:
\begin{equation}
	\begin{aligned}
		\inf &\; w \\
	    \text{s.t.} &\; f^r (\boldsymbol{p}) + d^{f}_{\boldsymbol{p}}(\boldsymbol{p})(\boldsymbol{y}-\boldsymbol{p}) \leq w, \forall \boldsymbol{p} \in \mathscr{Y} \\
	    &\; \boldsymbol{y} \in \mathscr{Y}.	
	\end{aligned}
\end{equation}
\end{proposition}

Thus, for the latest iterative solution $\boldsymbol{p}$, set the $SP(\boldsymbol{p})$ as the subproblem and  
$f^r (\boldsymbol{p}) + d^{f}_{\boldsymbol{p}}(\boldsymbol{p})(\boldsymbol{y}-\boldsymbol{p}) \leq w$ as the generated cut.
For each iteration, we add the latest iterative solution $\boldsymbol{p}$ into $\mathscr{U}$ 
which is the set of all results during the iterations.

Through the framework of Benders Decomposition method, we could obtain the optimal solution by iteratively solving the problem, 
denoted as $\mathbb{F}(\mathscr{U})$: 

\begin{align*}
	\inf &\; w \\
	\text{s.t.} &\; f^r (\boldsymbol{p}) + d^{f}_{\boldsymbol{p}}(\boldsymbol{p})(\boldsymbol{y}-\boldsymbol{p}) \leq w,, \forall \boldsymbol{p} \in \mathscr{U}\\
	&\; \boldsymbol{y} \in \mathscr{Y}.
\end{align*}

The detailed algorithm is given below.

\begin{algorithm}[H]
	\caption{RO}\label{algorithm1}
	Select any $\bar{\boldsymbol{y}} \in \mathscr{Y}$, solve $SP(\bar{\boldsymbol{y}})$ to obtain $f(\bar{\boldsymbol{y}})$, generate cuts, and add to $\mathbb{F}(\mathscr{U})$\;
	Set the $f_{best} \leftarrow f(\bar{\boldsymbol{y}})$,$\boldsymbol{y}_{best} \leftarrow \bar{\boldsymbol{y}}$, optimality gap $gap \leftarrow +inf$, the number of iteration $b \leftarrow 1$ and the set $\mathscr{U} = \{\bar{\boldsymbol{y}}\}$\;
	\While{($gap > \epsilon$)and($time < T$)}{Solve the problem $\mathbb{F}(\mathscr{U})$ with the optimal value $w^*$\;
	\If{a feasible vector, denoted by $\bar{\boldsymbol{y}}$ is not founded}{break}
	Solve $SP(\bar{\boldsymbol{y}})$ to obtain $f(\bar{\boldsymbol{y}})$ and calculate 
	$d^f_{\bar{\boldsymbol{y}}}(\bar{\boldsymbol{y}})$, generate cuts, and add to $\mathbb{F}(\mathscr{U})$\;
	Set $\mathscr{U} = \mathscr{U} \cup \{\bar{\boldsymbol{y}}\}, gap \leftarrow f(\bar{\boldsymbol{y}})-w^*$ and $b \leftarrow b+1$\;
	\If{$f_{best} \le f(\bar{\boldsymbol{y}})$}{$\boldsymbol{y}_{best} \leftarrow \bar{\boldsymbol{y}}$}
	}
    
\end{algorithm}

Especially, in the first step, we could use the optimal solution obtained from the deterministic graph 
with the respective means as arc lengths through the classical greedy algorithms like Prim algorithm.

With the solution framework shown above, as long as we could figure out $f(\boldsymbol{s})$ and $d^f_{\boldsymbol{s}}(\boldsymbol{s})$, 
we would obtain the optimal spanning tree in a finite number of steps. 
This conclusion is drawn from the fact that $\mathscr{U}$ is finite and $\mathscr{Y}$ would increase by one element each iteration.

\subsection{The calculation of the subgradient}

Here, we give out some important results and then the formulations of $f(\boldsymbol{s})$ and $d^f_{\boldsymbol{s}}(\boldsymbol{s})$.

\begin{proposition}\label{prop:1}
	The worst-case certainty equivalent has some useful properties as is shown below:
    \begin{description}
		\item[1.] \textbf{Monotonicity:} $C_{\alpha}(\tilde{t})$ is decreasing in $\alpha\ge 0$ and strictly decreasing when $\tilde{t}$ is not constant. Moreover,
		$$\lim_{\alpha \downarrow 0}C_{\alpha}(\tilde{t})=\overline{t_{\mathbb{F}}}, 
		\lim_{\alpha \to \inf}C_{\alpha}(\tilde{t})=\sup_{\mathbb{P}\in\mathbb{F}}\mathbb{E}_{\mathbb{P}}(\tilde{t}),$$
		where $\overline{t_{\mathbb{F}}}=\inf\{t\in\mathscr{R}|\mathbb{P}(\tilde{t}\le t)=1, \forall\mathbb{P}\in\mathbb{F}\}$; 
		\item[2.] \textbf{Convexity:} For any $\lambda\in[0,1], \tilde{t_1}, \tilde{t_2}\in\mathcal{V}$, and $\alpha_1, \alpha_2\ge 0$,
		$$C_{\lambda\alpha_1+(1-\lambda)\alpha_2}(\lambda\tilde{t_1}+(1-\lambda)\tilde{t_2})\le \lambda C_{\alpha_1}(\tilde{t_1}) + (1-\lambda)C_{\alpha_2}(\tilde{t_2});$$
		\item[3.] \textbf{Additivity:} If the random variables $\tilde{t_1}, \tilde{t_2}\in\mathcal{V}$ are independent of each other, then for any $\alpha\ge 0$,
		$$C_{\alpha}(\tilde{t_1}+\tilde{t_2})=C_{\alpha}(\tilde{t_1})+C_{\alpha}(\tilde{t_2}).$$
	\end{description}
\end{proposition}

Property (1) means that the smaller $\alpha$ is, the larger $C_{\alpha}(\cdot)$ will be. Property (2) shows that $C_{\alpha}(\tilde{t})$ is jointly convex in $(\alpha, \tilde{t})$. And property (3) indicates that $C_{\alpha}(\tilde{t})$ would be additive if random variables are independent from each other. The properties are attractive since we could derive the convexity of function $\rho$ directly from the property of $C_{\alpha}(\tilde{t})$, which is of great importance for the context of risk management where we wish that any convex combination of feasible solutions would be preferred too. 
To be specific, in our context where weight edge random variables are independent, we have that 
$C_{\alpha}(\tilde{\boldsymbol{w}}'\boldsymbol{y})=\sum_{e \in E}C_{\alpha}(\tilde{w}_ey_e)=\sum_{e \in E}C_{\alpha}(\tilde{w}_e)y_e$. The first equality is 
due to the additivity and the second is because that $\boldsymbol{y}$ is the incidence vector.
Plus, as we would analyze later, all of them would be beneficial to build a tractable model, according to which we could propose efficient algorithms.

\begin{rem}
	In the following, we would always consider the case with independent random variables without any note for simplicity, in spite of the fact that in practice, 
	the uncertainty would always be correlated. For example, uncertain travel times would be determined by some common factors, e.g. weather conditions, coincidence 
	of traffic jams. In the extension of \cite{jaillet2016routing}, they discuss the situation a little and propose a possible way to make extension. 
	The correlated situation is worth of more detailed discussion while here we only consider the independent case for simplicity to avoid the 
	tremendous increase in modeling and computational complexity. For example, \cite{qi2016preferences} has proven that 
	the path selection problem that Minimizes the certainty equivalent of total travel time
	$$\min_{\boldsymbol{s}\in\mathscr{S}_{SP}}C_{\alpha}(\boldsymbol{z}'\boldsymbol{s})$$ is NP-hard when the arc travel times are correlated.
\end{rem}

By the convexity and additivity property, the subproblem above can be a convex problem and solved efficiently because of its linear objective function. The direct approach is bisection method and other convex optimization algorithms could also be used. 

To calculate the subgradient $d^f_{\boldsymbol{s}}(\boldsymbol{s})$, we first write out its Lagrangian function
$$L(\boldsymbol{y},\alpha, \lambda) = \alpha + \lambda (h(\alpha, \boldsymbol{y})-\tau).$$

The constraints are in fact Slater's conditions, ensuring that the strong duality of the model holds. Then, we have for any $\boldsymbol{y} \in \mathscr{Y}$,

$$f^r (\boldsymbol{y})=\sup_{\lambda \ge 0}(\inf_{\alpha \ge 0} L(\boldsymbol{y},\alpha, \lambda))$$

\begin{proposition}\label{prop1}
	The subgradient of $f(\boldsymbol{p})$ with respect to $p_a$ for all $a \in \mathcal{A}$ can be calculated as 
	$$
	d^{f}_{p_a}(\boldsymbol{p})=d^{L}_{p_a}(\boldsymbol{p},\alpha^{*}, \lambda^{*})=
	\begin{cases}
	0, \quad & \alpha^* = 0, a \in \mathcal{A} \\
	-\frac{d^c_{p_a}(\alpha^*,\boldsymbol{p})}{d^c_{\alpha}(\alpha^*,\boldsymbol{p})}, \quad & \alpha^* >0,  a \in \mathcal{A} 
	\end{cases}
	$$
	Function $d^c_{p_a}(\alpha^*,\boldsymbol{p})$ and $d^c_{\alpha}(\alpha^*,\boldsymbol{p})$ is the subgradient of 
	$C_{\alpha}(\tilde{\boldsymbol{w}}'\boldsymbol{p})$ with respect to $p_a$ and $\alpha$ at point $(\alpha^*,\boldsymbol{p})$.
\end{proposition}

Next, we give out the detailed formulation of $f(\boldsymbol{s})$ and $d^f_{\boldsymbol{s}}(\boldsymbol{s})$ for specific uncertainty sets.

\begin{exmp}\label{uncertainty}
	If the distributional uncertainty set of random variable $\tilde{z}_a$ is a continuous random variable, 
	and the uncertainty set of random variable $\tilde{z}_a$  is given as below:
	$$ \mathbb{F}_a=\left\{
	\mathbb{P} \;\left|\; E_{\mathbb{P}}[\tilde{z}_a] \in [\underline{\mu}_a,\overline{\mu}_a],
	\mathbb{P}(\tilde{z}_a \in [\underline{z}_a,\overline{z}_a])=1
	\right.
	\right\}$$
	then
	$$
	C_{\alpha}(\tilde{z}_a) = \sup_{\mathbb{P} \in \mathbb{F}}\alpha ln E_{\mathbb{P}}(exp(\frac{\tilde{z}_a}{\alpha})) = 
	\begin{cases}
		\alpha ln(g(\tilde{z}_a)exp(\frac{\underline{z}_a}{\alpha})+h(\tilde{z}_a)exp(\frac{\overline{z}_a}{\alpha})), \quad & \alpha >0, \\
		\overline{z}_a, \quad & \alpha = 0.
	\end{cases}
	$$
	where $g(\tilde{z}_a)=\frac{\overline{z}_a-\overline{\mu}_a}{\overline{z}_a-\underline{z}_a}$ and 
	$h(\tilde{z}_a)=\frac{\overline{\mu}_a-\underline{z}_a}{\overline{z}_a-\underline{z}_a}$.
\end{exmp}

Immediately, as the function $C_{\alpha}(\tilde{\boldsymbol{w}}'\boldsymbol{p})$ is differentiable, we calculate its gradient with respect to $p_a$ as
$$
d^c_{p_a}(\alpha,\boldsymbol{p})=\frac{\partial}{\partial p_a}C_{\alpha}(\tilde{w}_a p_a)=
\frac{g(\tilde{z}_a)exp(\underline{z}_a p_a/{\alpha})\underline{z}_a+h(\tilde{z}_a)exp(\overline{z}_a p_a/{\alpha})\overline{z}_a}{g(\tilde{z}_a)exp(\underline{z}_a p_a/{\alpha})+h(\tilde{z}_a)exp(\overline{z}_a p_a/{\alpha})}.
$$

Meanwhile, the gradient of $C_{\alpha}(\tilde{\boldsymbol{w}}'\boldsymbol{p})$ with respect to $\alpha$ is 

\begin{align*}
	 d^c_{\alpha}(\alpha,\boldsymbol{p}) &\;=\sum_{a \in \mathcal{A}} \frac{\partial}{\partial \alpha}C_{\alpha}(\tilde{w}_a'p_a) \\
	 &\;= \sum_{a \in \mathcal{A}} (\ln(g(\tilde{z}_a)exp(\frac{\underline{z}_a}{\alpha})+h(\tilde{z}_a)exp(\frac{\overline{z}_a}{\alpha}))-\frac{g(\tilde{z}_a)exp(\underline{z}_a p_a/{\alpha})\underline{z}_a+h(\tilde{z}_a)exp(\overline{z}_a p_a/{\alpha})\overline{z}_a}{g(\tilde{z}_a)exp(\underline{z}_a p_a/{\alpha})+h(\tilde{z}_a)exp(\overline{z}_a p_a/{\alpha})} \frac{p_a}{\alpha}).
\end{align*}

\begin{exmp}\label{list}
For more kinds of ambiguity set, we list the equivalent representations of $C_\alpha(\lambda \tilde{z})$ below:
\begin{table}[H]
	\centering
	\scalebox{1}{
	    \begin{tabular}{c|c}
	    \hline
	    Ambiguity set&$C_\alpha(\lambda \tilde{z})=\sup_{P \in \mathbb{F}}\alpha\ln E_P[\exp(\frac{\lambda\tilde{z}}{\alpha})]$\\
	    \hline
	    $\mathbb{F}_1$ & $\min \{\alpha\log(\frac{(1+\underline{\mu})e^{\lambda/\alpha}+(1-\underline{\mu})e^{-\lambda/\alpha}}{2}), \alpha\log(\frac{(1+\overline{\mu})e^{\lambda/\alpha}+(1-\overline{\mu})e^{-\lambda/\alpha}}{2})\}$\\
		\hline
	    $\mathbb{F}_2$ &$ \alpha\log(\frac{\delta}{2(\mu+1)}e^{-\lambda/\alpha}+\frac{\delta}{2(1-\mu)}e^{\lambda/\alpha}+(1-\frac{\delta}{2(\mu+1)}-\frac{\delta}{2(1-\mu)})e^{\mu\lambda/\alpha})$\\
		\hline
		$\mathbb{F}_3$&
		$\min \{\alpha\log(\frac{{(1-\mu)}^2\exp(\frac{(\mu-\sigma^2)\lambda}{(1-\mu)\alpha})+(\sigma^2-\mu^2)\exp(\lambda/\alpha)}{1-2\mu+\sigma^2}),
		\alpha\log(\frac{{(1+\mu)}^2\exp(\frac{(\mu+\sigma^2)\lambda}{(1+\mu)\alpha})+(\sigma^2-\mu^2)\exp(-\lambda/\alpha)}{1+2\mu+\sigma^2})\}$\\
		\hline
	    $\dots$&$\dots$\\
	    \hline	
		\end{tabular}}
		\caption{\textbf{Equivalent representations of $C_\alpha(\lambda \tilde{z})$}}
\end{table}

where
$$\mathbb{F}_1 = \left\{
		\begin{array}{rcl}
		& E_{\mathbb{P}}[\tilde{\boldsymbol{w}}]\in [\underline{\boldsymbol{\mu}},\bar{\boldsymbol{\mu}}]\\	
		& \mathbb{P}[\tilde{\boldsymbol{w}}\in[-1,1]]=1\end{array}\right\}$$

$$\mathbb{F}_2 = \Biggl\{
		\begin{array}{rcl}
		& E_{\mathbb{P}}[\tilde{\boldsymbol{w}}]= \boldsymbol{\mu}\\
		& E_{\mathbb{P}}[|\tilde{\boldsymbol{w}}-\mu|]\le\delta\\
		& \mathbb{P}[\tilde{\boldsymbol{w}}\in[-1,1]]=1\end{array}\Biggr\}$$

$$\mathbb{F}_3 = \Biggl\{
	\begin{array}{rcl}
	& E_{\mathbb{P}}[\tilde{\boldsymbol{w}}]= \boldsymbol{\mu}\\
	& E_{\mathbb{P}}[|\tilde{\boldsymbol{w}}|^2]\le\sigma^2\\
	& \mathbb{P}[\tilde{\boldsymbol{w}}\in[-1,1]]=1\end{array}\Biggr\}$$
\end{exmp}

\begin{rem}
	Without loss of generality, we assume that the support is normalized. Otherwise, we could set $z \mapsto \frac{2z-(\underline{z}+\overline{z})}{\overline{z}-\underline{z}}$, 
	if $\mathbb{P}[\tilde{z}\in[\underline{z},\overline{z}]]=1$. For more interesting cases, readers could find them in \cite{chen2021robust}.
	And for the following contents, we mainly focus on the special case of the uncertainty set as in Example \ref{uncertainty}, 
	i.e. 
	$$ \mathbb{F}_a=\left\{
	\mathbb{P} \;\left|\; E_{\mathbb{P}}[\tilde{z}_a] = \mu,
	\mathbb{P}(\tilde{z}_a \in [\underline{z}_a,\overline{z}_a])=1
	\right.
	\right\}$$
	It's the special case when $\underline{\mu}=\overline{\mu}$.
\end{rem}

\subsection{Repeated-Prim(RP) algorithm} 

We find that by the framework of Benders Decomposition, the iterative convergence
is slow. That's because the calculation of $d^f_{\boldsymbol{y}}$ is unfortunately subgradient
of $f(\boldsymbol{y})$ rather than gradient, which means that the $f$ value of the new generated 
spanning tree isn't necessarily better than the ones before. In fact, it could be much worse
than even the corresponding value of average weight in our experiments. For example, it could happen that
$$\alpha^{(k+n)} \gg \alpha^{k}$$ for most of the $n$ where $\alpha^{k}$ is the risk parameter obtained from the configuration
through the k'th iteration while we wish that $\alpha^{k}$ would decrease continuingly and as fast as possible.

However, the RV Index enjoys perfect properties and the MST problem has quick algorithm like
Prim or Kruskal algorithm which could get results greatly faster than the formulation
of integer linear programming. Luckily, we could take advantage of these merits to 
generate exact results in greatly shorter time.

Especially, we know that the problem is equivalent to 
$$\inf_{\boldsymbol{s}\in\mathscr{Y}}f(\boldsymbol{s})$$
where when $\boldsymbol{s}$ is fixed, 
$$f(\boldsymbol{s}) = \{\alpha | \sum_{a}C_{\alpha}(\tilde{w}_a)s_a \le \tau\}.$$

Directly, for any given $\alpha\ge 0$ we could solve 
$\min_{\boldsymbol{y}\in\mathscr(Y)}C_{\alpha}(\tilde{\boldsymbol{w}}'\boldsymbol{y})$
by standard MST algorithms like Prim algorithm, and then we could use bisection algorithm 
to find the optimal $\alpha$, which is similar to the situation of shortest path problems based on minimizing the \textit{RV Index} as is pointed out in \cite{jaillet2016routing}. Thus we could get the conclusion below:

\begin{proposition}
	The RMST problem is polynomial solvable when the random variable $\tilde{\boldsymbol{w}}$ are independent of each other.
	More specifically, it could be solved in finite iterations of polynomial time algorithm, Prim algorithm.
\end{proposition}

In this case, we could use some tricks to improve the rate of convergence. Here, we replace the last weight by $C_{\alpha ^*}(\tilde{w}_a)$ where $\alpha ^*= f(\bar{\boldsymbol{s}})$ 
and $\bar{\boldsymbol{s}}$ is the spanning tree generated in 
last iteration by Prim algorithm. We could continue the procedure until the result could 
bring no improvement.

\begin{algorithm}[H]
	\caption{$RP$}\label{algorithm2}
	\KwData{graph$\leftarrow$\{edge:statistical information\}}
	\KwResult{Configuration $\boldsymbol{s}_{new}$}

	Initialization: edges$\leftarrow$\{edge:mean\},  $\boldsymbol{s}_{old} \leftarrow {None}$\;
	$\boldsymbol{s}_{new} \leftarrow$ $Prim(edges)$\;
	\While{ $\boldsymbol{s}_{new}\neq \boldsymbol{s}_{old}$}{
		$\alpha \leftarrow f(\boldsymbol{s}_{new})$\;
		calculate $C_{\alpha}(\tilde{w}_a)$ for all $a \in E$\;
		edges$\leftarrow$\{edge:$C_{\alpha}(\tilde{w}_a)$\}\;
		$\boldsymbol{s}_{old} \leftarrow \boldsymbol{s}_{new}$\;
		$\boldsymbol{s}_{new} \leftarrow$ $Prim(edges)$
	}
    
\end{algorithm}

\begin{proposition}\label{prop2}
Let $\boldsymbol{s}^k$ denotes the configuration obtained through the k step. We have that$f(\boldsymbol{s}^{k+1}) \le f(\boldsymbol{s}^k)$. 
\end{proposition}

\begin{proposition}\label{prop3}
	By the \textit{RP} algorithm, we could get the exact result $\boldsymbol{s}_{exact}$ in a finite number of steps. 
\end{proposition}

\section{Computational Study}

In this section, we perform experiments to evaluate whether
the RV Index model is practical solvable and whether the RV Index criterion can provide us a reasonable solution under uncertainty. 
First, we check the performance of the two proposed algorithms. Then, we use other benchmarks that minimizes average weight and maximizes budget of uncertainty
respectively and compare their performance with our model. Finally, we check the robustness against the changed parameters in the network and the target level. 
The program is coded in Python and run on a Intel Core i7 PC 
with a 2.21 GHz CPU by calling Gurobi as ILP solver.

\subsection{Performance of the algorithms}

Because the weight parameters were 
described by a range and a mean value, we assume that the weight distribution is piece-wise uniform. Instances with $10 \le |V| \le 30$ were used.
We carry out the first experiment to make a comparative study on the efficiency of the two algorithms proposed above. For a randomly generated 
Erdős–Rényi model $G(n, p)$, we solve the RV Index model by the two algorithms and compare the average CPU time and the times of iteration under 3 randomly generated Instances.  In order to ensure the problem feasibility, we artificially set the
target as 
$\tau = (1-\beta)\min_{\boldsymbol{s}\in\mathscr{Y}}{\boldsymbol{\mu}'\boldsymbol{s}}+\beta\min_{\boldsymbol{s}\in\mathscr{Y}}{\boldsymbol{\bar{\boldsymbol{z}}}'\boldsymbol{s}}$.
In this example, $\beta = 0.2$.

\begin{table}[!htbp]
\centering
\caption{Performances of Algorithms for the Target-based distributionally robust MST model.}
\begin{tabular}{c|c|c|c|c}
\hline

\multirow{3}*{$|V|$}&\multicolumn{4}{c}{$performance\  measures$}\\
\cline{2-5}
{ }&\multicolumn{2}{c|}{$RO$}&\multicolumn{2}{c}{$RP$}\\
\cline{2-5}
{ }&$CPU(s)$&$iteration$&$CPU(s)$&$iteration$\\
\hline
10&0.03&1.667&0.002&1.0\\
%\hline
20&20.932&40.333&0.009&1.667\\
%\hline
30&3728.522&268.667&0.014&2.0\\
\hline	
\end{tabular}	
\end{table}

Through the table above, we could know that compared with $RO$ algorithm, $RP$ algorithm could greatly shorten the CPU time needed to 
obtain results and the number of iteration while the $RO$ algorithm converges slowly
which is impracticable in the reality. However, things are different when the original case is more complex. Once in the deterministic case the 
efficient algorithm like Prim algorithm is unable to be used, the $RO$ algorithm could at least afford us basic practicable solution and enjoy great 
value of theoretical research. Thus in terms of the following comparative experiment, we would always take advantage of the $RP$ algorithm.

Next, we would compare our RP algorithm with direct bisection method. Because of the accessibility of the standard MST algorithm, now we could greatly expand our networks 
and enlarge the number of trial instances used in experiments. We randomly generate 50 instances and compare the statistics on CPU time of these two algorithms
for a network with 300 nodes. Table 4.2 suggests the calculation time of RP algorithm is much shorter than the bisection method. It provides an encouraging result for the 
employment of RP algorithm in the stochastic MST problem.

\begin{table}[H]
	\centering
	\caption{Statistics of CPU time of two algorithms.}
	\begin{tabular}{c|c|c|c}
	\hline
	
	\multirow{2}*{$Statistics$}&$Bisection$&\multicolumn{2}{c}{$RP$}\\
	\cline{2-4}
	{ }&$CPU\  time (sec)$&$CPU\  time (sec)$&$Number\  of\  iterations$\\
	\hline
	$Average$&3.7918&1.2874&3.84\\
	%\hline
	$Maximum$&6.9&2.5509&5.0\\
	%\hline
	$Minimum$&1.6376&0.4907&3.0\\
	$Standard\  deviation$&1.5751&0.509&0.5095\\
	\hline	
	\end{tabular}	
\end{table}

\subsection{Comparative study}

Next, we research on whether the RV Index criterion could bring us a reasonable solution under uncertainty and how $RP$ algorithm would 
perform with changed parameters.

\textbf{Benchmarks for stochastic minimum spanning tree problem}

We carry out the experiment to make a comparative study on the validity of the RV Index. For a randomly generated 
network, we solve a minimum spanning tree problem with target $\tau$ under uncertainty. We investigate several classical selection criteria
to find optimal paths. We summarize three criteria which appeared in the literature. And the setting is similar to the situation before.

\textbf{\emph{Minimize average weight}}

For a network with uncertain weight, the simplest way to find a spanning tree
is by minimizing the average weight, which can be formulated as 
a deterministic minimum spanning tree problem.

$$\min_{\boldsymbol{s} \in \mathscr{Y}} \boldsymbol{\mu}'\boldsymbol{s}.$$

\textbf{\emph{Maximize arrival probability}}

The second selection criterion is to find a path
that gives the largest probability to realize the target,
which is formulated as follows:
$$\max_{\boldsymbol{s}\in\mathscr{Y}}\mathbb{F}(\tilde{\boldsymbol{w}}'\boldsymbol{s}\leq\tau)$$

Since the problem is intractable, we adopt a sampling average approximation method
to solve it. Assuming the sample size is $K$, then we solve
\begin{align*}
	\max &\; \frac{1}{K}\sum_{k=1}^{K}I_k \\
	\text{s.t.} &\; \boldsymbol{s}'\boldsymbol{w}^k \le M(1-I_k)+\tau, &\,k=1,...,K, \\
	&\; I_k \in \{0,1\}, &\,k=1,...,K,\\
	&\; \boldsymbol{s} \in \mathscr{Y}.
\end{align*}

\begin{rem}
 This formulation has some insufferable demerits. Firstly, with respect to the choice of $M$, for smaller $M$, we could just exclude actual optimal solution while for larger $M$, it would take longer computational time. It has the same problem when it comes to the choice of $K$. What's worse, the randomly generated scenarios $\boldsymbol{w}^k$ would bring great difference to the final solutions, which is unsatisfactory.
\end{rem}

\textbf{\emph{Maximize budget of uncertainty}}

With the robust formulation proposed by \cite{bertsimas2003robust},
the robust minimum spanning tree problem is given as:
$$\min_{\boldsymbol{s}\in\mathscr{Y}}\max_{\tilde{\boldsymbol{w}}\in\mathscr{W}_{\Gamma}}\tilde{\boldsymbol{w}}'\boldsymbol{s}$$
in which, $\mathscr{W}_{\Gamma}=\{\boldsymbol{\mu}+\boldsymbol{c} |
\boldsymbol{0}\le\boldsymbol{c}\le\bar{\boldsymbol{w}}-\boldsymbol{\mu},
\sum_{a\in\mathscr{A}}\frac{c_a}{\bar{w}_a-\mu_{a}}\le\Gamma\}$, for all $\Gamma\ge0$.

Given the target $\tau$, we transform the problem to find a spanning tree that can return the maximal 
$\Gamma$ while respecting the target. The formulation is given as 
\begin{align*}
	\Gamma^* = \max &\; \Gamma \\
	\text{s.t.} &\;\max_{\tilde{\boldsymbol{w}}\in\mathscr{W}_{\Gamma}}\tilde{\boldsymbol{w}}'\boldsymbol{s} \le \tau, \\
	&\; \boldsymbol{s} \in \mathscr{Y}.
\end{align*}

By the mathematical framework suggested by \cite{bertsimas2003robust}, we get
$$\Gamma^* = \max_{l=1,\dots,|\mathscr{A}|+1}\frac{\tau - C_l}{\bar{w}_l-\mu_l}$$
where $C_l = \min_{\boldsymbol{s}\in\mathscr{Y}}(\boldsymbol{\mu}'\boldsymbol{s}+\sum_{j=1}^{l}((\bar{w}_j-\mu_j)-(\bar{w}_l-\mu_l))s_j),
l = 1,\dots,|\mathscr{A}|+1$ and $\mathscr{A}$ is the set of all edges in the network.

\begin{rem}
	Theoretically, the framework would be more complicated than our model, since we must execute the Prim algorithm for at least $|\mathscr{A}|+1$ or roughly $\frac{n(n-1)}{2}p$ times, which would accumulate sharply with the expansion of the networks. Here $(n, p)$ is the parameter of the Erdős–Rényi model $G(n, p)$.
\end{rem}
	
\begin{rem}
	 In this paper, we choose the similar criteria as the numerical experiences in \cite{jaillet2016routing}. 
	 %For later researches, we recommend to take more related criteria like the chance constraint model as is shown in our literature review.
\end{rem}

\textbf{Comparative study on the stochastic minimum spanning tree problem}

For maximizing the arrival probability,we just exclude it. That's because we use a sampling average approximation, which could lead to 
inconsistent solutions for comparison and sharp increase in computational complexity.  

For each instance, we randomly generate an undirected network with 300 nodes. The corresponding upper and lower bounds are randomly generated.
Table below summarizes the average performances among 50 instances. For notational clarity,
except the failure probability and CPU time, we only show the performance
ratio, which is the original performance divided by the performance of minimizing
the RV Index. Thus, the performance ratios for the RV Index model are one. The reason
that we would maintain the value of failure probability and CPU time is that we 
want to show the real situation of the performance rather than comparative situation.

\begin{table}[H]\label{table3}
	\centering
	\caption{Performances of various selection criteria for stochastic MST problem.}
	\scalebox{0.75}{
		\begin{threeparttable}
	    \begin{tabular}{c|c|c|c|c|c|c|c|c}
	    \hline
	    \multirow{2}*{$Selection\  criteria$}&\multicolumn{8}{c}{$performance\  measures$}\\
	    \cline{2-9}
	    { }&$Mean$&$Failure\  probability$&$STDEV$&$EL$&$CEL$&$VaR@95\%$&
	    $VaR@99\%$&$CPU\  time$\\
	    \hline
	    $Minimize\  average\  weight$&0.9933&0.04&1.7932&45.6869&2.7718&1.0076&1.0103&0.0914\\
	    %\hline
	    %Maximize arrival probability& & & & & & & & \\
	    %\hline
	    $Maximize\  budget\  of \ uncertainty$&0.9989&0.033&1.5084&44.6682&2.2356&1.0045&1.0097&1640.2159\\
	    %\hline
	    $Minimize\  the\  RV\  Index$&1.0&0.002&1.0&1.0&1.0&1.0&1.0&1.9333\\
	    \hline	
		\end{tabular}
			\begin{tablenotes}
				\footnotesize
				\item[1] \textbf{STDEV} refers to standard deviation;
				\item[2] \textbf{EL} refers to expected failure, $EL=E_P((\tilde{\boldsymbol{w}}'\boldsymbol{s}^*-\tau)^{+})$;
				\item[3] \textbf{CEL} refers to conditional expected failure, $CEL=E_P((\tilde{\boldsymbol{w}}'\boldsymbol{s}^*-\tau)^{+}|\tilde{\boldsymbol{w}}'\boldsymbol{s}^*\ge\tau)$;
				\item[4] \textbf{VaR@$\gamma$} refers to value-at-risk, $VaR@\gamma=\inf\{\nu\in\mathscr{R}|P(\tilde{\boldsymbol{w}}'\boldsymbol{s}^*>\nu)\le1-\gamma\}$. 
			\end{tablenotes} 
		\end{threeparttable}
	}
\end{table}

When compared with other selection criteria, the RV Index enjoys great improvement, especially on the failure probability. Although with higher 
mean value of the weight, the RV Index model enjoys more robustness and lower failure probability to meet the target, which is what we want. With
respect to the model of maximizing budget of uncertainty, it sacrifices the calculation time to reach robustness which is unacceptable to employ in 
the reality. Moreover, thanks to the greatly low failure probability, the \textbf{EL} value of minimizing
the RV Index is well below the level of the others, which could be more than 40 times than the prior. The \textbf{STDEV} and \textbf{CEL} value also mirrors the 
priority of minimizing the RV Index, being approximately the half of the others' levels.
Although the value of \textbf{VaR@$\gamma$} would mirror little difference, that could make sense given both the original
low failure probability and the risk allowed by the measure itself.

%\subsection{Robustness}

\subsubsection{Comparative study with the expansion of networks}

To highlight the superiority of the RV Index model, we change the number of nodes in the networks from 100 to 250 and we summarize the average 
performance among 10 instances. Again, we use the performance ratio to measure the performance of mean and standard deviation. Moreover, given
that it would cost a great amount of time to obtain optimal solution of the Maximizing budget of uncertainty when the networks are big,
we show the log formulation of the CPU time.

\begin{figure}[H]
	\centering
	\caption{Performances comparison with changing number of nodes.}
	\includegraphics[width=15cm]{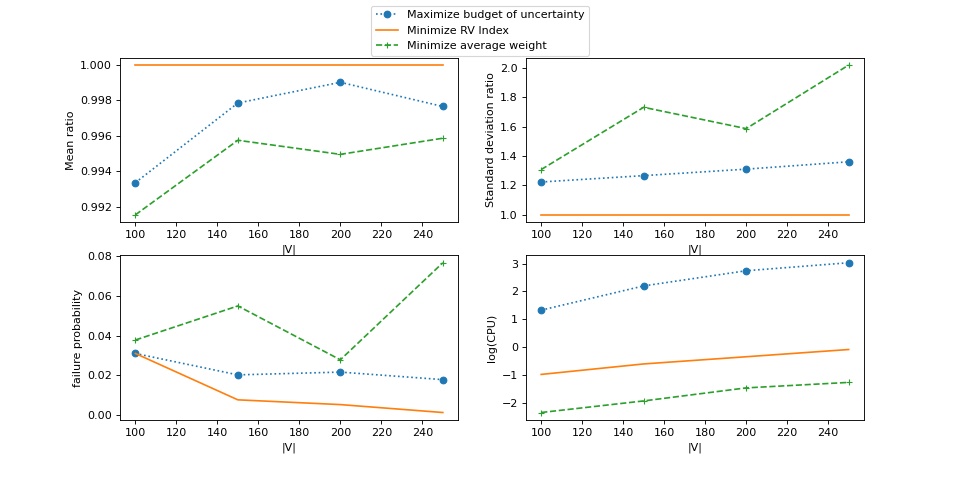}
\end{figure}

Apparently, the RV Index model is greatly superior to other models. With the expansion of the networks, the RV Index model would obtain
improved solution when faced with uncertainty in the weight of edges, which is reflected by the obviously lower level of failure probability
to meet the target and standard deviation. Although the maximizing budget of uncertainty enjoys more robustness compared with the minimizing 
average weight model, it cost too much calculation time which is unusually significant even in the log formulation of CPU time while the increase
of time corresponding to the expansion of networks of the RV Index model is steady and slow.

\subsubsection{Comparative study with changing target}

By varying the coefficient $\eta$, we also alter the target and summarize the performance ratio of each selection criterion.
Given that when $\eta \ge 0.3$, the failure probability would greatly converge to 0, we only show the result of failure probability ratio
and Expected lateness ratio with $\eta < 0.3$.

\begin{figure}[H]
	\centering
	\caption{Performances comparison with changing $\eta$.}
	\includegraphics[width=15cm]{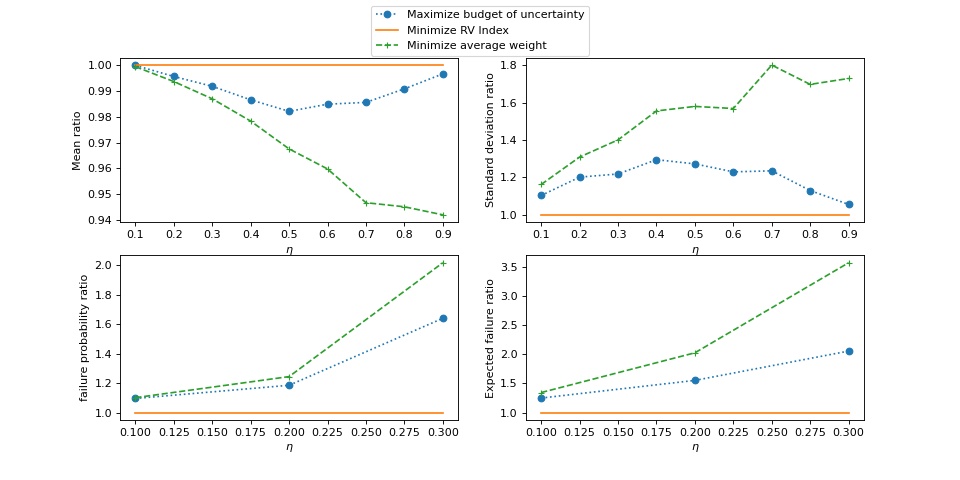}
\end{figure}

Among the remaining selection criteria, the RV Index model outperforms the others.
It's worthwhile to point out that it's not surprising that with increasing
$\eta$, the result of the Maximizing budget of uncertainty model would converge 
to the result of the Minimizing RV Index model in contradiction with the minimizing
average weight model. It's because of the given target, with which the minimizing average 
weight model would in fact have nothing to do while the result of the others would
be influenced by the set target.

\section{Conclusion}

In this paper, we first present a comprehensive summary of the literature related to the minimum spanning tree
problems, including its various deterministic variants, stochastic and robust spanning tree problems, and then propose a 
target-based distributionally robust optimization framework to deal with the minimum spanning tree problem in stochastic networks where the probability distribution function of the edge weight/cost is unknown but some statistical
information could be involved. We also propose two solution approaches based on the Benders decomposition framework and Prim algorithm
respectively. The first approach takes advantage of the mixed integer programming formulation, leading to the theoretical 
value for the further extension to other variants of minimum spanning tree problem in stochastic graphs. The second approach 
modifies the classical greedy algorithm for the deterministic minimum spanning tree problem, called Prim algorithm. Thus it enjoys more satisfactory
algorithmic aspect especially when faced with large-scale networks. In our computational experiments, the solutions obtained by the 
distributionally robust approach outperforms the others when faced with the uncertainty in edge weights.

\section{Appendix: Proofs}
\subsubsection{Proof of Proposition \ref{prop1}}
\begin{proof}
	{Proof:}

	By proposition 4 in \cite{jaillet2016routing}, We know that $d^{f}_{\boldsymbol{p}}(\boldsymbol{p})=d^{L}_{\boldsymbol{p}}(\boldsymbol{p},\alpha^{*}, \lambda^{*})$,
	where $$(\alpha^{*}, \lambda^{*}) \in Z(\boldsymbol{p}) = \left\{(\alpha^{0}, \lambda^{0})|L(\boldsymbol{p},\alpha^{0}, \lambda^{0})=\sup_{\lambda \ge 0}(\inf_{\alpha \ge 0} L(\boldsymbol{p},\alpha, \lambda)) \right \}.$$
	
	That's because :
	\begin{align*}
	    f^r(\boldsymbol{s})-f^r(\boldsymbol{p}) &\;  =\sup_{\lambda\ge 0}(\inf_{\lambda\ge 0}L(\boldsymbol{s},\alpha,\lambda))-\sup_{\lambda\ge 0}(\inf_{\lambda\ge 0}L(\boldsymbol{p},\alpha,\lambda))\\
	    &\; \ge \inf_{\lambda\ge 0}L(\boldsymbol{s},\alpha,\lambda^*)-\inf_{\lambda\ge 0}L(\boldsymbol{p},\alpha,\lambda^*)\\
	    &\; \ge d^L_{\boldsymbol{p}}(\boldsymbol{p},\alpha,\lambda^*)(\boldsymbol{s}-\boldsymbol{p}).
	\end{align*}

	Note that if $\alpha^* > 0$, 
	$$\frac{\partial}{\partial p_a}L(\boldsymbol{p},\alpha, \lambda) = \lambda^* d^c_{p_a}(\alpha^*,\boldsymbol{p})$$
	According to the generalized KKT Theorem, $\alpha ^{*}$ is primal optimal if and only if

	$$
	\begin{cases}
	1+\lambda^* d^c_{\alpha}(\alpha^*,\boldsymbol{p})=0, \quad &\\
	\lambda^*(h(\alpha^*, \boldsymbol{p})-\tau)=0, \quad &\\
	h(\alpha^*, \boldsymbol{p}) - \tau \leq 0, \quad &
	\end{cases}
	$$

	The equations above state that
	$$\lambda^* = - \frac{1}{d^c_{\alpha}(\alpha^*,\boldsymbol{p})}$$
	Hence we calculate the subgradient as 
	$$
	d^{f}_{p_a}(\boldsymbol{p})=d^{L}_{p_a}(\boldsymbol{p},\alpha^{*}, \lambda^{*})=
	\begin{cases}
	0, \quad & \alpha^* = 0, a \in \mathcal{A} \\
	-\frac{d^c_{p_a}(\alpha^*,\boldsymbol{p})}{d^c_{\alpha}(\alpha^*,\boldsymbol{p})}, \quad & \alpha^* >0,  a \in \mathcal{A} 
	\end{cases}
	$$
	
\end{proof}

\subsubsection{Proof of Example \ref{uncertainty} and \ref{list}}
\begin{proof}
	{Proof:} 
	For Example \ref{uncertainty}, we could conclude the proof by the following procedure.
	$$\sup_{\mathbb{P}\in\mathbb{F}}\alpha\ln E_{\mathbb{P}}(\exp(\frac{\tilde{z}}{\alpha}))=
	\alpha\ln\sup_{\mathbb{P}\in\mathbb{F}}E_{\mathbb{P}}(\exp(\frac{\tilde{z}}{\alpha}))=
	\alpha\ln(g(\tilde{z})\exp(\exp(\frac{\underline{z}}{\alpha}))+h(\tilde{z})\exp(\exp(\frac{\overline{z}}{\alpha}))).$$
	The optimal value is achieved when the random variable is subject to the two-point distribution, i.e.
	$$\mathbb{P}(\tilde{z}=\overline{z})=\frac{\mu-\underline{z}}{\overline{z}-\underline{z}}, 
	\mathbb{P}(\tilde{z}=\underline{z})=\frac{\overline{z}-\mu}{\overline{z}-\underline{z}}$$
	For Example \ref{list}, We only give out the brief proof for $\mathbb{F}_2$ and the others' could be similar.
\begin{align*}
	\sup_{\mathbb{P}\in\mathbb{F}}E_{\mathbb{P}}[e^{\lambda\tilde{z}}] \le&\;  \inf_{\gamma\ge0, \alpha, \beta} \alpha+\beta\mu+\gamma\delta\\
	\textit{s.t.}&\;  e^{\lambda z}\le\alpha+\beta z+\gamma |z-\mu|, \  \forall z \in [-1, 1]\\
	= &\; \inf_{\gamma\ge0, \alpha, \beta} \alpha+\beta\mu+\gamma\delta\\
	\textit{s.t.}&\;  e^{\lambda z}\le\alpha+\beta z+\gamma (z-\mu), \  \forall z \in [\mu, 1]\\
	&\;   e^{\lambda z}\le\alpha+\beta z+\gamma (\mu -z), \  \forall z \in [-1, \mu]\\
	=&\; \inf_{\gamma\ge0, \alpha, \beta} \alpha+\beta\mu+\gamma\delta\\
	\textit{s.t.}&\; e^{\lambda} \le \alpha+\beta+\gamma (1-\mu) \\
	&\; e^{\lambda\mu} \le \alpha+\beta \mu\\
	&\; e^{-\lambda} \le \alpha-\beta+\gamma (1+\mu) \\
	= &\; \sup_{p_1, p_2, p_3\ge 0} p_1e^{\lambda}+p_2e^{\lambda\mu}+p_3e^{-\lambda}\\
	\textit{s.t.}&\; p_1+p_2+p_3=1\\
	&\; p_1+p_2\mu-p_3=\mu\\
	&\; (1-\mu)p_1+(1+\mu)p_3 \le \delta
\end{align*}
The first inequality is due to weak duality and the second equality is by the optimal solution of a 
convex maximization problem attained at the boundary. The third equality is because of linear optimization strong duality.
By solving the last linear optimization problem, we could draw the conclusion.
\end{proof}

\subsubsection{Proof of Proposition \ref{prop2}}
\begin{proof}
	{Proof:}If $\boldsymbol{s}^{k+1} = \boldsymbol{s}^{k}$, it's trivial.

Otherwise, let $\alpha^{k}$ is optimal for $(\boldsymbol{s}^k, \tau)$, and we know $\sum_{a}C_{\alpha^k}(\tilde{w}_a)s^k_a = \tau$ by the definition.

Plus $\sum_{a}C_{\alpha^k}(\tilde{w}_a)s^{k+1}_a \le \sum_{a}C_{\alpha^k}(\tilde{w}_a)s^k_a$ because of the generation procedure of $\boldsymbol{s}^{k+1}$.

If $\sum_{a}C_{\alpha^k}(\tilde{w}_a)s^{k+1}_a < \sum_{a}C_{\alpha^k}(\tilde{w}_a)s^k_a$,
which means that $\sum_{a}C_{\alpha^k}(\tilde{w}_a)s^{k+1}_a < \tau$.

By the fact that
$C_{\alpha}(\cdot)$ is decreasing and continuous in $\alpha \ge 0$, we must have 
$\alpha^{k+1} = f(\boldsymbol{s}^{k+1}) < \alpha^k = f(\boldsymbol{s}^k)$.

If $\sum_{a}C_{\alpha^k}(\tilde{w}_a)s^{k+1}_a = \sum_{a}C_{\alpha^k}(\tilde{w}_a)s^k_a = \tau$,
we naturally have that $f(\boldsymbol{s}_{k+1}) = f(\boldsymbol{s}^k)$.
\end{proof}

\subsubsection{Proof of Proposition \ref{prop3}}
\begin{proof}
	{Proof:}Assume that for fixed $\boldsymbol{s}$(the last generated one), we have optimal $\alpha^* = f(\boldsymbol{s})$.
So now we replace the weight of edge $a$ by $C_{\alpha ^*}(\tilde{w}_a)$. By Prim algorithm,
now we could get another spanning tree $\boldsymbol{s}_{next}$. By the proposition above, we only need to prove that if $\boldsymbol{s}_{next} = \boldsymbol{s}$, it must be
the exact result $arg \min_{\boldsymbol{s} \in \mathscr{Y}}f(\boldsymbol{s})$.

We only need to prove that 
$\sum_{a}C_{\alpha^*_{next}}(\tilde{w}_a)s^{next}_a \ge \sum_{a}C_{\alpha^*_{next}}(\tilde{w}_a)s^{exact}_a$.

If $\boldsymbol{s}_{next} \neq \boldsymbol{s}_{exact}$, we have 
$\alpha^*_{exact} < \alpha^*_{next}$, which means that 
$\sum_{a}C_{\alpha^*_{next}}(\tilde{w}_a)s^{exact}_a < \sum_{a}C_{\alpha^*_{exact}}(\tilde{w}_a)s^{exact}_a = \tau$.

Thus we have $\sum_{a}C_{\alpha^*_{next}}(\tilde{w}_a)s^{exact}_a < \tau = \sum_{a}C_{\alpha^*_{next}}(\tilde{w}_a)s^{next}_a$.
It contradicts.

We could get the exact solution in a finite number of steps because the optimal $\alpha$ would be decreased strictly in each iteration. Once the optimal $\alpha$ stops to be improved, the solution would have been obtained by the discussion above. And we know that the whole solution space is limited and each optimal $\alpha$ would correspond to one spanning tree in the solution space.
\end{proof}

\newpage

\bibliographystyle{plainnat}
\bibliography{References}

\begin{thebibliography}{50}
\providecommand{\natexlab}[1]{#1}
\providecommand{\url}[1]{\texttt{#1}}
\expandafter\ifx\csname urlstyle\endcsname\relax
  \providecommand{\doi}[1]{doi: #1}\else
  \providecommand{\doi}{doi: \begingroup \urlstyle{rm}\Url}\fi

\bibitem[Aggarwal et~al.(1982)Aggarwal, Aneja, and Nair]{aggarwal1982minimal}
Vijay Aggarwal, Yash~P Aneja, and KPK Nair.
\newblock Minimal spanning tree subject to a side constraint.
\newblock \emph{Computers \& Operations Research}, 9\penalty0 (4):\penalty0 287--296, 1982.

\bibitem[Akbari~Torkestani and Meybodi(2012)]{akbari2012learning}
Javad Akbari~Torkestani and Mohammad~Reza Meybodi.
\newblock A learning automata-based heuristic algorithm for solving the minimum spanning tree problem in stochastic graphs.
\newblock \emph{The Journal of Supercomputing}, 59\penalty0 (2):\penalty0 1035--1054, 2012.

\bibitem[Alexopoulos and Jacobson(2000)]{alexopoulos2000state}
Christos Alexopoulos and Jay~A Jacobson.
\newblock State space partition algorithms for stochastic systems with applications to minimum spanning trees.
\newblock \emph{Networks: An International Journal}, 35\penalty0 (2):\penalty0 118--138, 2000.

\bibitem[Aron and Van~Hentenryck(2004)]{aron2004complexity}
Ionu{\c{t}}~D Aron and Pascal Van~Hentenryck.
\newblock On the complexity of the robust spanning tree problem with interval data.
\newblock \emph{Operations Research Letters}, 32\penalty0 (1):\penalty0 36--40, 2004.

\bibitem[Bertsimas and Sim(2003)]{bertsimas2003robust}
Dimitris Bertsimas and Melvyn Sim.
\newblock Robust discrete optimization and network flows.
\newblock \emph{Mathematical programming}, 98\penalty0 (1):\penalty0 49--71, 2003.

\bibitem[Bertsimas et~al.(2011)Bertsimas, Brown, and Caramanis]{bertsimas2011theory}
Dimitris Bertsimas, David~B Brown, and Constantine Caramanis.
\newblock Theory and applications of robust optimization.
\newblock \emph{SIAM review}, 53\penalty0 (3):\penalty0 464--501, 2011.

\bibitem[Bui and Zrncic(2006)]{bui2006ant}
Thang~N Bui and Catherine~M Zrncic.
\newblock An ant-based algorithm for finding degree-constrained minimum spanning tree.
\newblock In \emph{Proceedings of the 8th annual conference on Genetic and evolutionary computation}, pages 11--18, 2006.

\bibitem[Chen and Sim(2021)]{chen2021robust}
Li~Chen and Melvyn Sim.
\newblock Robust cara optimization.
\newblock \emph{Available at SSRN 3937474}, 2021.

\bibitem[Chen and Sim(2009{\natexlab{a}})]{chen2009goal}
Wenqing Chen and Melvyn Sim.
\newblock Goal-driven optimization.
\newblock \emph{Operations Research}, 57\penalty0 (2):\penalty0 342--357, 2009{\natexlab{a}}.

\bibitem[Chen and Sim(2009{\natexlab{b}})]{Chen2009}
W.Q. Chen and M.~Sim.
\newblock Goal-driven optimization.
\newblock \emph{Operations Research}, 57\penalty0 (2):\penalty0 342--357, 2009{\natexlab{b}}.

\bibitem[Chen et~al.(2009)Chen, Hu, and Hu]{chen2009polynomial}
Xujin Chen, Jie Hu, and Xiaodong Hu.
\newblock A polynomial solvable minimum risk spanning tree problem with interval data.
\newblock \emph{European Journal of Operational Research}, 198\penalty0 (1):\penalty0 43--46, 2009.

\bibitem[Cui et~al.(2021)Cui, Long, QI, and Zhang]{cui2021inventory}
Zheng Cui, Daniel~Zhuoyu Long, Jin QI, and Lianmin Zhang.
\newblock Inventory routing problem under uncertainty.
\newblock \emph{Available at SSRN 3946807}, 2021.

\bibitem[de~Almeida et~al.(2005)de~Almeida, Yamakami, and Takahashi]{de2005evolutionary}
T~Agostinho de~Almeida, Akebo Yamakami, and M{\'a}rcia~Tomie Takahashi.
\newblock An evolutionary approach to solve minimum spanning tree problem with fuzzy parameters.
\newblock In \emph{International Conference on Computational Intelligence for Modelling, Control and Automation and International Conference on Intelligent Agents, Web Technologies and Internet Commerce (CIMCA-IAWTIC'06)}, volume~2, pages 203--208. IEEE, 2005.

\bibitem[Doulliez and Jamoulle(1972)]{doulliez1972transportation}
P~Doulliez and E~Jamoulle.
\newblock Transportation networks with random arc capacities.
\newblock \emph{Revue fran{\c{c}}aise d'automatique, informatique, recherche op{\'e}rationnelle. Recherche op{\'e}rationnelle}, 6\penalty0 (V3):\penalty0 45--59, 1972.

\bibitem[Frieze and Tkocz(2021)]{frieze2021randomly}
Alan~M Frieze and Tomasz Tkocz.
\newblock A randomly weighted minimum arborescence with a random cost constraint.
\newblock \emph{Mathematics of Operations Research}, 2021.

\bibitem[Gamvros et~al.(2006)Gamvros, Golden, and Raghavan]{doi:10.1287/ijoc.1040.0123}
Ioannis Gamvros, Bruce Golden, and S.~Raghavan.
\newblock The multilevel capacitated minimum spanning tree problem.
\newblock \emph{INFORMS Journal on Computing}, 18\penalty0 (3):\penalty0 348--365, 2006.
\newblock \doi{10.1287/ijoc.1040.0123}.
\newblock URL \url{https://doi.org/10.1287/ijoc.1040.0123}.

\bibitem[Gouveia et~al.(2011)Gouveia, Simonetti, and Uchoa]{gouveia2011modeling}
Luis Gouveia, Luidi Simonetti, and Eduardo Uchoa.
\newblock Modeling hop-constrained and diameter-constrained minimum spanning tree problems as steiner tree problems over layered graphs.
\newblock \emph{Mathematical Programming}, 128\penalty0 (1):\penalty0 123--148, 2011.

\bibitem[Gruber et~al.(2006)Gruber, van Hemert, and Raidl]{gruber2006neighbourhood}
Martin Gruber, Jano van Hemert, and G{\"u}nther~R Raidl.
\newblock Neighbourhood searches for the bounded diameter minimum spanning tree problem embedded in a vns, ea, and aco.
\newblock In \emph{Proceedings of the 8th annual conference on Genetic and Evolutionary Computation}, pages 1187--1194, 2006.

\bibitem[Hall et~al.(2015)Hall, Long, Qi, and Sim]{Hall2015}
N.G. Hall, D.Z. Long, J.~Qi, and M.~Sim.
\newblock Managing underperformance risk in project portfolio selection.
\newblock \emph{Operations Research}, 63\penalty0 (3):\penalty0 660--675, 2015.

\bibitem[He and Qi(2008)]{he2008model}
Fangguo He and Huan Qi.
\newblock A model and algorithm for minimum spanning tree problems in uncertain networks.
\newblock In \emph{2008 3rd International Conference on Innovative Computing Information and Control}, pages 493--493. IEEE, 2008.

\bibitem[Hochbaum(2003)]{hochbaum2003efficient}
Dorit~S Hochbaum.
\newblock Efficient algorithms for the inverse spanning-tree problem.
\newblock \emph{Operations Research}, 51\penalty0 (5):\penalty0 785--797, 2003.

\bibitem[Hutson and Shier(2005)]{hutson2005bounding}
Kevin~R Hutson and Douglas~R Shier.
\newblock Bounding distributions for the weight of a minimum spanning tree in stochastic networks.
\newblock \emph{Operations research}, 53\penalty0 (5):\penalty0 879--886, 2005.

\bibitem[Hutson and Shier(2006)]{hutson2006minimum}
Kevin~R Hutson and Douglas~R Shier.
\newblock Minimum spanning trees in networks with varying edge weights.
\newblock \emph{Annals of Operations Research}, 146\penalty0 (1):\penalty0 3--18, 2006.

\bibitem[Ishii and Matsutomi(1995)]{ishii1995confidence}
H~Ishii and T~Matsutomi.
\newblock Confidence regional method of stochastic spanning tree problem.
\newblock \emph{Mathematical and computer modelling}, 22\penalty0 (10-12):\penalty0 77--82, 1995.

\bibitem[Ishii and Nishida(1983)]{ishii1983stochastic}
Hiroaki Ishii and Toshio Nishida.
\newblock Stochastic bottleneck spanning tree problem.
\newblock \emph{Networks}, 13\penalty0 (3):\penalty0 443--449, 1983.

\bibitem[Ishii et~al.(1981)Ishii, Shiode, Nishida, and Namasuya]{ishii1981stochastic}
Hiroaki Ishii, Sh{\=o}go Shiode, Toshio Nishida, and Yoshikazu Namasuya.
\newblock Stochastic spanning tree problem.
\newblock \emph{Discrete Applied Mathematics}, 3\penalty0 (4):\penalty0 263--273, 1981.

\bibitem[Jaillet et~al.(2016)Jaillet, Qi, and Sim]{jaillet2016routing}
Patrick Jaillet, Jin Qi, and Melvyn Sim.
\newblock Routing optimization under uncertainty.
\newblock \emph{Operations research}, 64\penalty0 (1):\penalty0 186--200, 2016.

\bibitem[Jain and Mamer(1988)]{jain1988approximations}
Anjani Jain and John~W Mamer.
\newblock Approximations for the random minimal spanning tree with application to network provisioning.
\newblock \emph{Operations Research}, 36\penalty0 (4):\penalty0 575--584, 1988.

\bibitem[Kruskal(1956)]{10.2307/2033241}
Joseph~B. Kruskal.
\newblock On the shortest spanning subtree of a graph and the traveling salesman problem.
\newblock \emph{Proceedings of the American Mathematical Society}, 7\penalty0 (1):\penalty0 48--50, 1956.
\newblock ISSN 00029939, 10886826.
\newblock URL \url{http://www.jstor.org/stable/2033241}.

\bibitem[Li and Deshpande(2011)]{li2011maximizing}
Jian Li and Amol Deshpande.
\newblock Maximizing expected utility for stochastic combinatorial optimization problems.
\newblock In \emph{2011 IEEE 52nd Annual Symposium on Foundations of Computer Science}, pages 797--806. IEEE, 2011.

\bibitem[Long and Qi(2014)]{long2014distributionally}
Daniel~Zhuoyu Long and Jin Qi.
\newblock Distributionally robust discrete optimization with entropic value-at-risk.
\newblock \emph{Operations Research Letters}, 42\penalty0 (8):\penalty0 532--538, 2014.

\bibitem[Magnanti and Wolsey(1995)]{magnanti1995optimal}
Thomas~L Magnanti and Laurence~A Wolsey.
\newblock Optimal trees.
\newblock \emph{Handbooks in operations research and management science}, 7:\penalty0 503--615, 1995.

\bibitem[Mao(1970)]{mao1970survey}
James~CT Mao.
\newblock Survey of capital budgeting: Theory and practice.
\newblock \emph{Journal of finance}, pages 349--360, 1970.

\bibitem[Mohd(1994)]{mohd1994interval}
Ismail~Bin Mohd.
\newblock Interval elimination method for stochastic spanning tree problem.
\newblock \emph{Applied Mathematics and Computation}, 66\penalty0 (2-3):\penalty0 325--341, 1994.

\bibitem[Montemanni(2006)]{montemanni2006benders}
Roberto Montemanni.
\newblock A benders decomposition approach for the robust spanning tree problem with interval data.
\newblock \emph{European Journal of Operational Research}, 174\penalty0 (3):\penalty0 1479--1490, 2006.

\bibitem[Montemanni and Gambardella(2005)]{montemanni2005branch}
Roberto Montemanni and Luca~Maria Gambardella.
\newblock A branch and bound algorithm for the robust spanning tree problem with interval data.
\newblock \emph{European Journal of Operational Research}, 161\penalty0 (3):\penalty0 771--779, 2005.

\bibitem[{\"O}ncan(2007)]{oncan2007design}
Temel {\"O}ncan.
\newblock Design of capacitated minimum spanning tree with uncertain cost and demand parameters.
\newblock \emph{Information Sciences}, 177\penalty0 (20):\penalty0 4354--4367, 2007.

\bibitem[{\"O}ncan et~al.(2008){\"O}ncan, Cordeau, and Laporte]{oncan2008tabu}
Temel {\"O}ncan, Jean-Francois Cordeau, and Gilbert Laporte.
\newblock A tabu search heuristic for the generalized minimum spanning tree problem.
\newblock \emph{European Journal of Operational Research}, 191\penalty0 (2):\penalty0 306--319, 2008.

\bibitem[Parsa et~al.(1998)Parsa, Zhu, and Garcia-Luna-Aceves]{parsa1998iterative}
Mehrdad Parsa, Qing Zhu, and JJ~Garcia-Luna-Aceves.
\newblock An iterative algorithm for delay-constrained minimum-cost multicasting.
\newblock \emph{IEEE/ACM transactions on networking}, 6\penalty0 (4):\penalty0 461--474, 1998.

\bibitem[Paul et~al.(2020)Paul, Freund, Ferber, Shmoys, and Williamson]{paul2020budgeted}
Alice Paul, Daniel Freund, Aaron Ferber, David~B Shmoys, and David~P Williamson.
\newblock Budgeted prize-collecting traveling salesman and minimum spanning tree problems.
\newblock \emph{Mathematics of Operations Research}, 45\penalty0 (2):\penalty0 576--590, 2020.

\bibitem[Prim(1957)]{prim1957shortest}
Robert~Clay Prim.
\newblock Shortest connection networks and some generalizations.
\newblock \emph{The Bell System Technical Journal}, 36\penalty0 (6):\penalty0 1389--1401, 1957.

\bibitem[Qi et~al.(2016)Qi, Sim, Sun, and Yuan]{qi2016preferences}
Jin Qi, Melvyn Sim, Defeng Sun, and Xiaoming Yuan.
\newblock Preferences for travel time under risk and ambiguity: Implications in path selection and network equilibrium.
\newblock \emph{Transportation Research Part B: Methodological}, 94:\penalty0 264--284, 2016.

\bibitem[Rahimian and Mehrotra(2019)]{rahimian2019distributionally}
Hamed Rahimian and Sanjay Mehrotra.
\newblock Distributionally robust optimization: A review.
\newblock \emph{arXiv preprint arXiv:1908.05659}, 2019.

\bibitem[Salazar-Neumann(2007)]{salazar2007robust}
Martha Salazar-Neumann.
\newblock The robust minimum spanning tree problem: compact and convex uncertainty.
\newblock \emph{Operations research letters}, 35\penalty0 (1):\penalty0 17--22, 2007.

\bibitem[Shen et~al.(2015)Shen, Kurt, and Wang]{shen2015chance}
Siqian Shen, Murat Kurt, and Jue Wang.
\newblock Chance-constrained programming models and approximations for general stochastic bottleneck spanning tree problems.
\newblock \emph{INFORMS Journal on Computing}, 27\penalty0 (2):\penalty0 301--316, 2015.

\bibitem[Sokkalingam et~al.(1999)Sokkalingam, Ahuja, and Orlin]{sokkalingam1999solving}
PT~Sokkalingam, Ravindra~K Ahuja, and James~B Orlin.
\newblock Solving inverse spanning tree problems through network flow techniques.
\newblock \emph{Operations Research}, 47\penalty0 (2):\penalty0 291--298, 1999.

\bibitem[Sourd and Spanjaard(2008)]{sourd2008multiobjective}
Francis Sourd and Olivier Spanjaard.
\newblock A multiobjective branch-and-bound framework: Application to the biobjective spanning tree problem.
\newblock \emph{INFORMS Journal on Computing}, 20\penalty0 (3):\penalty0 472--484, 2008.

\bibitem[Torkestani and Meybodi(2011)]{torkestani2011learning}
Javad~Akbari Torkestani and Mohammad~Reza Meybodi.
\newblock Learning automata-based algorithms for solving stochastic minimum spanning tree problem.
\newblock \emph{Applied Soft Computing}, 11\penalty0 (6):\penalty0 4064--4077, 2011.

\bibitem[Wei et~al.(2021)Wei, Walteros, and Pajouh]{wei2021integer}
Ningji Wei, Jose~L Walteros, and Foad~Mahdavi Pajouh.
\newblock Integer programming formulations for minimum spanning tree interdiction.
\newblock \emph{INFORMS Journal on Computing}, 33\penalty0 (4):\penalty0 1461--1480, 2021.

\bibitem[Yaman et~al.(2001)Yaman, Kara{\c{s}}an, and P{\i}nar]{yaman2001robust}
Hande Yaman, Oya~Ekin Kara{\c{s}}an, and Mustafa~{\c{C}} P{\i}nar.
\newblock The robust spanning tree problem with interval data.
\newblock \emph{Operations research letters}, 29\penalty0 (1):\penalty0 31--40, 2001.

\end{thebibliography}

\end{document}